\documentclass[12pt]{iopart}
\usepackage{amsmath}
\usepackage{amssymb}
\usepackage{amsthm}
\usepackage{graphicx}
\def\N{\mathbb{N}}
\def\R{\mathbb{R}}
\def\xdag{x^\dagger}
\def\xd{x^\delta}
\def\yd{y^\delta}
\def\oh{\frac12}

\def\kp{k+}

\def\pp{\mathrm{p}}
\def\rr{\mathrm{r}}
\def\ppd{{\mathrm{p}^*}}
\def\rrd{{\mathrm{r}^*}}
\def\ppp{p}
\def\chf{1\!{\rm I}}
\usepackage{color}

\def\altil{\tilde{\alpha}}
\def\gamtil{\tilde{\gamma}}
\def\vtil{\tilde{v}}

\newtheorem{theorem}{Theorem}

\newtheorem{lemma}{Lemma}
\newtheorem{remark}{Remark}

\begin{document}
\title[Halley's method in Banach space]{A convergence rates result for an iteratively regularized Gauss-Newton-Halley method in Banach space
\thanks{}}
\author{B.~Kaltenbacher}
\address{Alpen-Adria Universit\"at Klagenfurt,\\
Universit\"atsstra\ss e 65--67, 9020 Klagenfurt, Austria\\
barbara.kaltenbacher@aau.at}

\begin{abstract}
The use of second order information on the forward operator often comes at a very moderate additional computational price in the context of parameter identification probems for differential equation models. On the other hand the use of general (non-Hilbert) Banach spaces has recently found much interest due to its usefulness in many applications. 
This motivates us to extend the second order method from \cite{HalleyNM}, (see also \cite{HettlichRundell}) to a Banach space setting and analyze its convergence. We here show rates results for a particular source condition and different exponents in the formulation of Tikhonov regularization in each step. 
This includes a complementary result on the (first order) iteratively regularized Gauss-Newton method (IRGNM) in case of a one-homogeneous data misfit term, which corresponds to exact penalization.
The results clearly show the possible advantages of using second order information, which get most pronounced in this exact penalization case.
Numerical simulations for a coefficient identification problem in an elliptic PDE illustrate the theoretical findings.  
\end{abstract}

\section{Introduction}
Identification of parameters in ordinary or partial differential equations by Newton methods usually requires repeated solution of the model equation (the PDE or ODE) and its linearization, since these methods rely on a first order Taylor expansion of the parameter-to-state map. It has already been observed in \cite{HettlichRundell,HalleyNM} that also higher derivative evaluation for this forward operator typically lead to the same linear differential equation as the one arising for the first order derivative, and only the right hand sides differ. Let us illustrate this by means of two examples.

\paragraph{Example 1}
Consider identification of the pair $(a,c)$ of possibly spatially varying coefficients in the nonlinear elliptic boundary value problem
\begin{equation}\label{ellPDE}
\left\{\begin{array}{rcl}
- \nabla f(a,\nabla u) + g(c,u)&=&0\mbox{ in }\Omega\\
u&=&h\mbox{ on }\partial\Omega
\end{array}\right.
\end{equation}
from measurements $y=Cu$ of the state $u$, where $C$ is some linear operator (e.g., a trace operator in case of boundary measurements).
Problems of this kind arise, e.g., in stationary inverse groundwater filtration, as well as in the characterization or nondestructive inspection of (non)linearly elastic or magnetic materials.
Here $\Omega\subseteq\R^d$, and the functions $f:\R^{d+1}\to\R$, $g:\R^2\to\R$, $h\in H^{1/2}(\partial\Omega)$ are given. Note that linear growth of $f,g$ and monotonicity (uniform one in case of $f$) with respect to their second argument allow to show well-posedness of \eqref{ellPDE} by means of the Lax-Milgram Lemma. Then, using the parameter-to-state map $G:(a,c)\mapsto u$, the forward operator and its derivatives at some point $(a,c)$ in parameter space can be written as $F(a,c)=C G(a,c)$, $F'(a,c)=CG'(a,c)$, $F''(a,c)=CG''(a,c)$, where the derivatives of $G$ at $(a,c)$ in certain directions can be recovered as solutions of the same linearized elliptic PDE
with different right hand sides: For parameter increments $(\alpha,\gamma)$, $(\altil,\gamtil)$ we have that $v^1=G'(a,c)(\alpha,\gamma)$, $v^2=G''(a,c)((\alpha,\gamma),(\altil,\gamtil))$ solve
\begin{equation}\label{ellPDElin}
\left\{\begin{array}{rcl}
- \nabla \Bigl(\partial_2 f(a,\nabla u) \, \nabla v^i\Bigr) + \partial_2 g(c,u) \, v^i
&=& b^i
 \mbox{ in }\Omega\\
v^i&=&0\mbox{ on }\partial\Omega
\end{array}\right.
\quad i=1,2\,,
\end{equation}
where  
\[
\begin{aligned}
b^1=&
\nabla \Bigl(\partial_1 f(a,\nabla u) \alpha \Bigr)- \partial_1 g(c,u) \gamma
\\[1ex]
b^2=
&-\partial_1^2 g(c, u) (\gamma,\gamtil) 
-\partial_1\partial_2 g(c, u) (\gamma,\vtil^1) 
-\partial_1\partial_2 g(c, u) (\gamtil, v^1) 
- \partial_2^2 g(c, u) (v^1,\vtil^1)\\
&+\nabla \Bigl(
\partial_1^2 f(a,\nabla u) (\alpha,\altil) 
+\partial_1\partial_2 f(a,\nabla u) (\alpha,\nabla\vtil^1) 
+\partial_1\partial_2 f(a,\nabla u) (\altil,\nabla v^1)\\
&\hspace*{9cm} 
+ \partial_2^2 f(a,\nabla u) (\nabla v^1,\nabla\vtil^1)
\Bigr)
\end{aligned}
\]
and $u=G(a,c)$, $\vtil^1=G'(a,c)(\altil,\gamtil)$,
i.e., the same linear elliptic boundary value problem \eqref{ellPDElin}, only with different right hand sides.

\paragraph{Example 2}
For modelling time dependent problems, consider the state space model
\begin{equation}\label{IVP}
\dot{u}(t)+ f(t,u(t),c) =0\,, \ t>0\,,\quad u(0)=u_0
\end{equation}
where the dot denotes the time derivative, which includes systems of ODEs but also (thinking of $u(t)$ as an element of a function space over some spatial domain $\Omega$) time dependent PDEs.
Given $f$, $u_0$, we seek to identify the paramter $c$ -- possibly element of a finite or infinite dimensional Banach space --
from measurements $y=Cu$ of the state $u$, where $C$ is some linear operator.
Again, using the parameter-to-state-map $G:c\mapsto u$, the forward operator and its derivatives at some point $c$ in parameter space can be written as $F(c)=C G(c)$, $F'(c)=CG'(c)$, $F''(c)=CG''(c)$, where $v^1=G'(c)\gamma$, $v^2=G''(c)(\gamma,\gamtil)$ solve the same linear system
\begin{equation}\label{IVPlin}
\dot{v}^i(t)+ \partial_2 f(t,u(t),c) v^i(t)=b^i\,, \ t>0\,,\quad v^i(0)=0
\quad i=1,2\,,
\end{equation}
with different right hand sides
\[
\begin{aligned}
b^1=&-\partial_3 f(t,u(t),c)\gamma
\\[1ex]
b^2=&
-\partial_3^2 f(t,u(t),c)(\gamma,\gamtil)
-\partial_2\partial_3 f(t,u(t),c)(\gamma,\vtil^1(t))
-\partial_2\partial_3 f(t,u(t),c)(\gamtil,v^1(t))\\
&-\partial_2^2 f(t,u(t),c)(v^1(t),\vtil^1(t))
\end{aligned}
\]
where $u=G(c)$, $\vtil^1=G'(c)\gamtil$.

\medskip

We point out that in these two examples (and many more), evaluating $F''(c_k)$ at some iterate $c_k$ leads to the same linear problem as evaluating $F'(c_k)$, just with some different right hand side. In case of elliptic or parabolic PDEs this means that the stiffness matrix remains unchanged and therefore, once we have done the computations for $F'(c_k)$, the additional effort for evaluating $F''(c_k)$ can be kept quite moderate, usually much lower than the effort for doing an additional Newton step that requires $F'(c_{k+1})$, thus setting up a new stiffness matrix, at some different coefficient $c_{k+1}$.
This cheap evaluation of the second derivative {\em at the same iterate} is just what Halley's method (see \cite{Brown77,Doering70,RenArgyros12}
for the well-posed setting) does, which for ill-posed problems in Hilbert spaces can be formulated as follows (the coefficient iterates are now denoted by $\xd_k$ instead of $c_k$):
\begin{equation}\label{HalleyIRGNM}
\begin{aligned}
&
x_0^\delta=x_0
\\
&
\text{for }k=1,2,\ldots
\\
&\hspace*{0.5cm}T_k= F'(\xd_k); \quad r_k= F(\xd_k)-\yd\\
&\hspace*{0.5cm}\xd_{\kp}=\xd_k-(T_k^*T_k+\beta_k I)^{-1}\{T_k^*r_k+\beta_k(\xd_k-x_0)\}\\
&\hspace*{0.5cm}S_k=T_k +\oh F''(\xd_k)(\xd_{\kp}-\xd_k,\cdot)\\
&\hspace*{0.5cm}\xd_{k+1}=\xd_k-(S_k^*S_k+\alpha_k I)^{-1}\{S_k^*r_k+\alpha_k(\xd_k-x_0)\}
\end{aligned}
\end{equation}
with two a priori fixed sequences $(\alpha_k)_{k\in\N}$, $(\beta_k)_{k\in\N}$ satisfying
\begin{equation}\label{q}
\alpha_k\searrow 0\,, \quad
\beta_k\searrow 0\,, \quad
1\leq \frac{\alpha_k}{\alpha_{k+1}}\leq q\,, \quad
1\leq \frac{\beta_k}{\beta_{k+1}}\leq q\,,
\end{equation}
cf. \cite{HalleyNM}.
In here, $F:X\to Y$ is the forward operator in the operator equation formulation 
\begin{equation}\label{Fxy}
F(x)=y
\end{equation}
of the coefficient identification problem and the superscript $\delta$ indicates the presence of noise in the given data $\yd$, whose deterministic level we assume to be known, i.e.,
\begin{equation}\label{delta}
\|y-\yd\|\leq\delta\,.
\end{equation}
Here we use a fixed reference point $x_0$, i.e., the first order version of this method (skipping the step with $S_k$) would be the iteratively regularized Gauss-Newton method (IRGNM), see, e.g., \cite{Baku92,HohageWerner13,IRGNMBanachrates}. If we would replace $x_0$ by the current iterate $\xd_k$ in each step, we would arrive at the Levenberg-Marquardt type version of Halley's method considered by Hettlich and Rundell in \cite{HettlichRundell}.
While \cite{HettlichRundell,HalleyNM} concentrate on the case of $X,Y$ being Hilbert spaces, it is often desirable to also work in Banach spaces that to not possess Hilbert space structure, such as $L^1$ or the space of Radon measures for obtaining sparse solutions or modelling impulsive noise, or $L^ \infty$ for guaranteeing essential bounds (e.g., nonnegativity) of coefficients or modelling uniform noise, cf., e.g., \cite{Clason, ClasonJin, HohageWerner14}.
Method \eqref{HalleyIRGNM} can indeed be extended to the more general setting of $X,Y$ being Banach spaces in a straightforward manner:
\begin{equation}\label{HalleyIRGNMBanach}
\begin{aligned}
&
x_0^\delta=x_0
\\
&
\text{for }k=1,2,\ldots
\\
&\hspace*{0.5cm}T_k= F'(\xd_k); \quad r_k= F(\xd_k)-\yd\\
&\hspace*{0.5cm}\xd_{\kp}\in\text{argmin}_x 
\frac{1}{\rr}\|T_k(x-\xd_k)+r_k\|^\rr+\frac{\beta_k}{\pp}\|x-x_0\|^\pp
\\
&\hspace*{0.5cm}S_k=T_k +\oh F''(\xd_k)(\xd_{\kp}-\xd_k,\cdot)\\
&\hspace*{0.5cm}\xd_{k+1}\in\text{argmin}_x 
\frac{1}{\rr}\|S_k(x-\xd_k)+r_k\|^\rr+\frac{\alpha_k}{\pp}\|x-\xd_k\|^\pp
\end{aligned}
\end{equation}
with $\pp,\rr\in[1,\infty)$.
Below we will prove a convergence result under the source condition 
\begin{equation}\label{scBanach}
T^*v\in J_\pp(\xdag-x_0)
\end{equation}
for some $v\in X$
Here $J_\pp=\partial \frac{1}{\pp}\|\cdot\|^\pp$ denotes the duality mapping.
Our analysis will make use of the shifted Bregman distance 
\[
D^{x_0}_{\pp,\xi}(\tilde{x},x)=\frac{1}{\pp}\|x-x_0\|^\pp-\frac{1}{\pp}\|\tilde{x}-x_0\|^\pp -\langle \xi,x-\tilde{x}\rangle \mbox{ with }\xi\in J_\pp(\tilde{x}-x_0)\,,
\]
which, if $X$ is $\pp$-convex, satisfies the coercivity estimate 
\begin{equation}\label{Xpcoercive} 
D^{x_0}_{\pp,\xi}(\tilde{x},x) \ge \underline c \,
\|\tilde x-x \|^\pp  \quad \mbox{for all} \quad \tilde x,x \in X 
\end{equation}
for some constant $\underline c>0$ depending on $\pp$ (see, e.g., \cite[Lemma~2.7]{Bonetal08}).
The stopping index $k_*$ will be the first one such that $\alpha_{k}^{\frac{1}{\rr-1}}\leq\tau\delta$, i.e., 
\begin{equation}\label{kstBanach}
\alpha_{k_*}^{\frac{1}{\rr-1}}\leq\tau\delta< \alpha_{k}^{\frac{1}{\rr-1}} 
\quad \forall k\in\{0,\ldots,k_*-1\}\,,
\end{equation}
if $\rr>1$. 
In case $\rr=1$ we can choose $\beta_k$ and $\alpha_k$ constant (i.e., the Tikhonov regularization parameter is independent of $\delta$, as typical for this case, see, e.g., \cite{BurgerOsher})
and 
\begin{equation}\label{kstBanach_r1}
k_*\geq\left[\log_\sigma(\log_2(\delta^{-1/\pp})))\right]\,.
\end{equation}
where $\sigma=\frac{\pp+2}{\pp^2}$, which is larger than one for $\pp<2$.

We first of all consider the case $\rr>1$.
\begin{theorem}\label{th:LipschitzBanachrgt1}
Let $X,Y$ be Banach spaces with additionally $X$ being $p$-convex so that \eqref{Xpcoercive} holds, where the exponents satisfy
\[
 1\leq \pp\leq 2 \mbox{ and } 1< \rr \leq \pp < 2\rr \,.
\]
 Assume that a source condition \eqref{scBanach} with  $\|v\|$ sufficiently small holds.
Let $F$ be twice Fr\'{e}chet differentiable with $F''$ bounded 
and Lipschitz continuous 
in $\mathcal{B}_\rho(\xdag)$ and let $x_0$ be suffiently close to $\xdag$. 
Assume that 
\begin{equation}\label{betaksalphak}
\beta_k=s\alpha_k 
\end{equation}
with $s>0$, $\alpha_0$ sufficiently small and \eqref{q}, 
and let $k_*$ be chosen according to \eqref{kstBanach} with $\tau$ sufficiently large.

Then the iterates defined by \eqref{HalleyIRGNMBanach} converge at the rate 
\[
\|\xd_{k_*}-\xdag\|=O(\delta^{\frac{1}{\pp}}) \mbox{ as }\delta\to0\,.
\]
If $\delta=0$ we have convergence
\begin{equation}\label{ratealpha2}
\|\xd_k-\xdag\|=O(\alpha_k^{\frac{1}{\pp(\rr-1)}}) \mbox{ as }k\to\infty\,.
\end{equation}
\end{theorem}

The case $\rr=1$ corresponding to exact penalization of the data misfit (see, e.g., \cite{BurgerOsher}) is treated separately. Since the existing results on the IRGNM from \cite{HohageWerner13,IRGNMBanachrates} do not seem to be applicable, \footnote{note that \cite{IRGNMBanachrates} requires $\rr>1$, whereas in \cite{HohageWerner13}, the multiplicative source condition (12) does not contain the case \eqref{scBanach} if $\rr=1$ and the additive one (35) with Theorem 4.2 does not give the desired rate in this case} we also prove the corresponding result for the IRGNM.

\begin{theorem}\label{th:LipschitzBanachreq1}
Let $X,Y$ be Banach spaces with additionally $X$ being $p$-convex so that \eqref{Xpcoercive} holds, where the exponents satisfy
\[
 1\leq \pp < 2 \mbox{ and }\rr=1\,.
\]
 Assume that a source condition \eqref{scBanach} with  $\|v\|$ sufficiently small holds.
Let $F$ be twice Fr\'{e}chet differentiable with $F''$ bounded 
and Lipschitz continuous 
in $\mathcal{B}_\rho(\xdag)$ and let $x_0$ be suffiently close to $\xdag$. 
Assume that 
\[
\alpha_k\geq\underline{\alpha}\,, \quad \beta_k\geq\underline{\beta}
\]
for some constants $\underline{\alpha},\underline{\beta}>0$, and let $k_*$ be chosen according to \eqref{kstBanach_r1} with $\sigma=\frac{\pp+2}{\pp^2}>1$.

Then the iterates defined by \eqref{HalleyIRGNMBanach} converge at the rate 
\begin{equation}\label{ratedeltaBanach_r1}
\|\xd_{k_*}-\xdag\|=O(\delta^{\frac{1}{\pp}}) \mbox{ as }\delta\to0\,.
\end{equation}
If $\delta=0$ we have convergence of order $\sigma=\frac{\pp+2}{\pp^2}$ (i.e., cubic for $\pp=1$). 

\smallskip

For the iteratively regularized Gauss-Newton method (IRGNM) defined by setting $x_{k+1}=x_{\kp}$, 
under the same assumptions (except the ones on $F''$ and on $\beta_k$, which are not needed) with $k_*$ chosen according to \eqref{kstBanach_r1} with $\sigma=\frac{2}{\pp}>1$, 
the iterates  converge at the rate \eqref{ratedeltaBanach_r1}. 
If $\delta=0$, we have convergence of order $\sigma=\frac{2}{\pp}$ for the IRGNM iterates (i.e., quadratic for $\pp=1$). 
\end{theorem}

\begin{remark}\label{rem:rates}
The difference between first and second order IRGNM becomes even clearer here than in the Hilbert space case with quadratic penalties from \cite{HalleyNM}, especially in the case $\rr=1$ of Theorem \ref{th:LipschitzBanachreq1}: 
In the exact data case, the order of convergence is always better for Halley than for IRGNM, since 
\[
\forall \pp\in[1,2) \, : \quad \sigma^{\mbox{\footnotesize Halley}}=\frac{\pp+2}{\pp^2}>\frac{2}{\pp}=\sigma^{\mbox{\footnotesize IRGNM}}\,,
\] 
which becomes most obvious in the case $\pp=1$, where we get cubic convergence for Halley's method and quadratic one for the IRGNM.
This faster convergence is also reflected in the number of iterates according to \eqref{kstBanach_r1} in case of noisy data, since the logarithm in \eqref{kstBanach_r1} is taken with respect to a larger basis for Halley than for IRGNM.
\end{remark}

\subsection{Proof of 
Theorems \ref{th:LipschitzBanachrgt1}, \ref{th:LipschitzBanachreq1}
}\label{subsecLipschitzBanach}
By minimality in the definition \eqref{HalleyIRGNMBanach} of the iterates we have 
\begin{equation}\label{estBanach1}
\frac{1}{\rr}\|K(x-\xd_k)+r_k\|^\rr+\frac{\kappa}{\pp}\|x-x_0\|^\pp
\leq\frac{1}{\rr}\|K(\xdag-\xd_k)+r_k\|^\rr+\frac{\kappa}{\pp}\|\xdag-x_0\|^\pp
\end{equation}
for 
\[
(x,K,\kappa)\in\{(\xd_{\kp},T_k,\beta_k)\,,\ (\xd_{k+1},S_k,\alpha_k)\}\,.
\]
On the other hand, we use the fact that by definition of the Bregman distance and the source condition \eqref{scBanach} in both cases we have
\begin{equation}\label{estBanach2}
\begin{aligned}
\frac{\kappa}{\pp}\|x-x_0\|^\pp-\frac{\kappa}{\pp}\|\xdag-x_0\|^\pp
&=\kappa D^{x_0}_{\pp,T^*v}(\xdag,x)+\kappa \langle T^*v,x-\xdag\rangle\\
&\geq \kappa D^{x_0}_{\pp,T^*v}(\xdag,x) - \kappa \|v\|\, \|T(x-\xdag)\|\,.
\end{aligned}
\end{equation}
Note that (up to a linearization error) we have
\begin{equation} \label{estBanach3}
\|K(x-\xd_k)+r_k\|^\rr \approx (\delta+\|T(x-\xdag)\|)^\rr \,, \quad \|K(\xdag-\xd_k)+r_k\|^\rr \approx \delta^\rr
\end{equation}
so we will correspondingly dominate the term $\kappa \|v\|\, \|T(x-\xdag)\|$ from \eqref{estBanach2} by a small multiple of $\|T(x-\xdag)\|^\rr$, which is obvious in case $\rr=1$ with the smallness assumption
\[
\kappa \|v\| < \frac{2^{-\rr}}{\rr}=\frac12
\]
and for $\rr\in(1,\infty)$ follows from Young's inequality in the form 
\begin{equation}\label{abeps}
a b\leq \epsilon a^\rr+ C(\epsilon,\rr)b^{\rrd}
\end{equation}
with 
\[
\rrd=\frac{\rr}{\rr-1}\,, \quad C(\epsilon,\rr)=\frac{\rr-1}{\epsilon^{1/(\rr-1)} \rr^{\rrd}}
\]
setting $\epsilon=2^{-\rr}/\rr$, $a=\|T(x-\xdag)\|$, $b=\kappa\|v\|$.
Putting these estimates together and using the simple inequalities 
$(a-b)^\rr+b^\rr\geq 2^{-(\rr-1)} a^\rr$, $(a+b)^\rr\leq 2^{\rr-1}(a^\rr+b^\rr)$ we obtain two estimates of the form
\[
\begin{aligned}
&2^{-(\rr-1)} \|T(x-\xdag)\|^\rr + \kappa D^{x_0}_{\pp,T^*v}(\xdag,x)\\
&\leq \kappa\|v\| \, \|T(x-\xdag)\| + 2(\delta + \text{Taylor remainder})^\rr\\
&\leq c^0 \|T(x-\xdag)\|^\rr + C^0 \|v\|^{\rrd} \kappa^{\rrd}  + (\delta + \text{Taylor remainder})^\rr
\end{aligned}
\]
where $0<c^0<2^{-(\rr-1)}$ and the term $\|v\|^{\rrd}\kappa^{\rrd} $ vanishes (also formally, by $\kappa\|v\|<1$) in case $\rr=1$, hence 
\[
\begin{aligned}
c^1 \|T(x-\xdag)\|^\rr + \kappa D^{x_0}_{\pp,T^*v}(\xdag,x)
\leq& C^1 \|v\|^{\rrd} \kappa^{\rrd} + 2(\delta + \text{Taylor remainder})^\rr
\end{aligned}
\]
for some constants $c^1, C^1 >0$ depending only on $q$.
The approximations in \eqref{estBanach3} can be quantified by the Taylor remainder estimates
\begin{equation*}
\begin{aligned}
&\|T_k(\xdag-\xd_k)+r_k\|\leq
\delta + \frac12 C_2\|\xd_k-\xdag\|^2\\[1ex]
&\|T_k(\xd_{\kp}-\xd_k)+r_k \ - \ T(\xd_{\kp}-\xdag)\|\leq 
\|T_k(\xdag-\xd_k)+r_k\|+\|t_1'\|\, \|\xd_{\kp}-\xdag\|\\
&\leq \delta + \frac12 C_2\|\xd_k-\xdag\|^2 + C_2\|\xd_k-\xdag\|\, \|\xd_{\kp}-\xdag\|\\[1ex]
&\|S_k(\xdag\xd_k)+r_k\|\leq 
\delta + \frac16 L_2\|\xd_k-\xdag\|^3 +\frac12 C_2 \|\xd_{\kp}-\xdag\|\,\|\xd_k-\xdag\|\\[1ex]
&\|S_k(\xd_{k+1}-\xd_k)+r_k \ - \ T(\xd_{k+1}-\xdag)\|\leq 
\|S_k(\xdag\xd_k)+r_k\|+\|t_2'\|\, \|\xd_{k+1}-\xdag\|\\
&\leq \delta + \frac16 L_2\|\xd_k-\xdag\|^3 +\frac12 C_2 \|\xd_{\kp}-\xdag\|\,\|\xd_k-\xdag\|\\
&\quad+\Bigl(\oh L_2\|\xd_k-\xdag\|^2+\oh C_2 (\|\xd_k-\xdag\|+\|\xd_{\kp}-\xdag\|)\Bigr)
\|\xd_{k+1}-\xdag\|\,.
\end{aligned}
\end{equation*}
Hence we end up with
\begin{equation}\label{estBanach4}
\begin{aligned}
&c\|T(\xd_{\kp}-\xdag)\|^\rr + \beta_k D^{x_0}_{\pp,T^*v}(\xdag,\xd_{\kp})
\\
&\leq C \Bigl( \|v\|^{\rrd} \beta_k^{\rrd}
+ \delta^\rr + \|\xd_k-\xdag\|^{2\rr} + \|\xd_k-\xdag\|^\rr\, \|\xd_{\kp}-\xdag\|^\rr
\Bigr)
\end{aligned}
\end{equation}
and
\begin{equation}\label{estBanach5}
\begin{aligned}
&c \|T(\xd_{k+1}-\xdag)\|^\rr + \alpha_k D^{x_0}_{\pp,T^*v}(\xdag,\xd_{k+1})\\
&\leq C \Bigl( \|v\|^{\rrd} \alpha_k^{\rrd}
+ \delta^\rr + \|\xd_k-\xdag\|^{3\rr} + \|\xd_{\kp}-\xdag\|^\rr\,\|\xd_k-\xdag\|^\rr\\
&\quad+\Bigl(\|\xd_k-\xdag\|^{2\rr}+\|\xd_k-\xdag\|^\rr+\|\xd_{\kp}-\xdag\|^\rr\Bigr)
\|\xd_{k+1}-\xdag\|^\rr
\Bigr)
\end{aligned}
\end{equation} 
with some constants $c,C>0$ depending only on $q,C_2,L_2$.
Thus, using the coercivity estimate \eqref{Xpcoercive} and considering first of all the case $\rr>1$, we expected to obtain the rates  
\[
\|\xd_{k+1}-\xdag\|=O(\alpha_k^{\frac{1}{\pp(\rr-1)}})\,, \quad 
\|\xd_{\kp}-\xdag\|=O(\beta_k^{\frac{1}{\pp(\rr-1)}})
\]
and hence consider the quantities
\begin{equation}\label{gammakGammakBanach}
\begin{aligned}
&\gamma_k=\frac{\|\xd_k-\xdag\|}{\alpha_k^{\frac{1}{\pp(\rr-1)}}}\,, \ k\in\{0,\ldots,k_*\}\,,\\
&
\Gamma_{k+1}=\frac{\|\xd_{\kp}-\xdag\|}{\beta_k^{\frac{1}{\pp(\rr-1)}}}\,, \ k\in\{1,\ldots,k_*-1\}\,, \ \Gamma_0=0\,,
\end{aligned}
\end{equation}
for which, dividing by $\beta_k^\rrd$ and $\alpha_{k+1}^\rrd$, respectively, from \eqref{estBanach4}, \eqref{estBanach5}, \eqref{Xpcoercive} and the stopping rule 
\begin{equation}\label{kstBanach_alphabeta}
\alpha_{k_*}\leq(\tau\delta)^{\rr-1}< \alpha_{k}\,, \
\beta_{k_*}\leq(\bar{\tau}\delta)^{\rr-1}< \beta_{k}
\quad \forall k\in\{0,\ldots,k_*-1\}\,,
\end{equation}
(cf. \eqref{kstBanach}) 
we obtain, for $k\leq k_*-1$
\begin{equation}\label{estBanach5a}
\begin{aligned}
\Gamma_{k+1}^\pp\leq \frac{C}{\underline{c}}\Bigl( \|v\|^{\rrd} + \frac{1}{\tau^\rr}
+ \underbrace{\alpha_k^{\frac{2\rrd}{\pp}}\beta_k^{-\rrd} \gamma_k^{2\rr}}_{I}
+ \underbrace{\alpha_k^{\frac{\rrd}{\pp}}\beta_k^{-\frac{\rrd}{\ppd}} \gamma_k^\rr\Gamma_{k+1}^\rr}_{II}
\Bigr)
\end{aligned}
\end{equation}
and
\begin{equation}\label{estBanach6}
\begin{aligned}
\gamma_{k+1}^\pp\leq& \frac{C}{\underline{c}}q^\rrd\Bigl( \|v\|^{\rrd} + \frac{1}{\tau^\rr}
+ \underbrace{\alpha_k^\frac{\rrd(3-\pp)}{\pp} \gamma_k^{3\rr}}_{III}
+ \underbrace{\beta_k^{\frac{\rrd}{\pp}}\alpha_k^{-\frac{\rrd}{\ppd}} \gamma_k^\rr\Gamma_{k+1}^\rr}_{IV}
\Bigr)\\
&+\frac{C}{\underline{c}}q^{\frac{\rrd}{\ppd}}\Bigl(
	\underbrace{\alpha_k^\frac{\rrd(3-\pp)}{\pp} \gamma_k^{2\rr}}_{V}
	+\underbrace{\alpha_k^\frac{\rrd(2-\pp)}{\pp} \gamma_k^{\rr}}_{VI}
	+\underbrace{\beta_k^{\frac{\rrd}{\pp}}\alpha_k^{-\frac{\rrd}{\ppd}} \Gamma_{k+1}^\rr}_{VII}
\Bigr) \gamma_{k+1}^{\rr}\,.
\end{aligned}
\end{equation}
Since the sequences $\alpha_k$, $\beta_k$ tend to zero, the desired boundedness of the right hand side imposes some restrictions to the exponents. Namely, in view of term VI in \eqref{estBanach6} we need
\begin{equation}\label{ple2}
\pp\leq 2
\end{equation}
and from terms I, II, IV, using the fact that by \eqref{ple2} $\frac{1}{\pp-1}\geq\frac{2}{\pp}$ we infer condition 
\begin{equation}\label{alphabetaBanach}
m \alpha_k^{\frac{2}{\pp}} \leq \beta_k \leq M \alpha_k^{\frac{1}{\pp-1}}
\end{equation}
for some $m,M >0$ independent of $k$. For instance, $\beta_k=s\alpha_k$ with some $s>0$ is an admissible choice satisfying \eqref{alphabetaBanach}, and setting $\bar{\tau}=s^{\frac{1}{\rr-1}}\tau$ in \eqref{kstBanach_alphabeta} guarantees well-definedness of $k_*$.
On the other hand, conditions \eqref{ple2}, \eqref{alphabetaBanach} imply boundedness of all the terms I-VII.
Therewith we end up with estimates
\begin{equation}\label{estBanach7}
\begin{aligned}
\Gamma_{k+1}^\pp\leq a + b \gamma_k^{2\rr} + c \gamma_k^\rr\Gamma_{k+1}^\rr
=: \phi(\gamma_k,\Gamma_{k+1})
\end{aligned}
\end{equation}
and
\begin{equation}\label{estBanach8}
\begin{aligned}
\gamma_{k+1}^\pp\leq& d+ e \gamma_k^{3\rr}+ f \gamma_k^\rr\Gamma_{k+1}^\rr
+ \Bigl(h \gamma_k^{2\rr}
+ i \gamma_k^{\rr}
+ j \Gamma_{k+1}^\rr
\Bigr) \gamma_{k+1}^{\rr}=:\Phi(\gamma_k,\gamma_{k+1},\Gamma_{k+1})
\,,
\end{aligned}
\end{equation}
where we have imposed the bound
\begin{equation}\label{d}
\|v\|^{\rrd}+\frac{1}{\bar{\tau}^\rr}\leq \frac{\underline{c}}{C} \min\{a, d q^{-\rrd}\}\,.
\end{equation}
We wish to carry out an induction proof of the claim 
\begin{equation}\label{induction1}
\gamma_k\leq\bar{\gamma}\,, \quad \Gamma_k\leq\bar{\Gamma} \quad  \forall k\in\{0,\ldots,k_*\}\,,
\end{equation}
(for all $k\in\mathbb{N}_0$ in case $\delta=0$) with appropriately chosen constants $\bar{\gamma}, \bar{\Gamma}>0$. For this purpose, it suffices to do the induction step, since the induction beginning can be easily established by imposing the closeness condition $\|\xd_0-\xdag\|\leq \bar{\gamma}\alpha_0^{\frac{1}{\pp(\rr-1)}}$ and using the convention $\Gamma_0=0$.
The induction step can be carried out by means of the following Lemma
\begin{lemma}\label{lem1}
Let $\phi$, $\Phi$ be defined as in \eqref{estBanach7}, \eqref{estBanach8} with 
\[
\rr\leq\pp\leq 2\rr\,.
\]
Then there exist $\bar{\gamma}$, $\bar{\Gamma}$ such that for $d,e,f,h,j>0$ sufficiently small the implication 
\[ 
\forall \gamma,\Gamma >0 \ : \quad 
\Bigl(\Gamma^\pp\leq \phi(\bar{\gamma},\Gamma)\mbox{ and }
\gamma^\pp\leq \Phi(\bar{\gamma},\gamma,\Gamma)\Bigr) \ \Rightarrow \
\Bigl(\Gamma\leq\bar{\Gamma}\mbox{ and }\gamma\leq\bar{\gamma}\Bigr)
\]
holds.
\end{lemma}
\begin{proof}
see the Appendix.
\end{proof}
Note that indeed $d$ can be made small by imposing $\|v\|$ small and choosing $\tau$ large;
$e,f,h,j$ can be made small by choosing 
$\beta_k=s \alpha_k$ with $s$ small, i.e., \eqref{betaksalphak}. 
Also note that the choice of all bounds in this lemma is independent of $\delta$ and of $k$.

Thus we have established \eqref{induction1}. This immediately implies the claimed rate in the exact data case.
The stopping rule \eqref{kstBanach} then implies the rate $\|\xd_{k_*}-\xdag\|=O(\delta^{\frac{1}{\pp}})$ in case of noisy data.

\medskip 

Finally, we consider the special case $\rr=1$ where \eqref{estBanach4}, \eqref{estBanach5} with the coercivity estimate \eqref{Xpcoercive} becomes
\begin{equation}\label{estBanach4_r1}
\begin{aligned}
&c\|T(\xd_{\kp}-\xdag)\| + \beta_k \underline{c}\|\xd_{\kp}-\xdag\|^\pp\\
&\leq C \Bigl(\delta + \|\xd_k-\xdag\|^{2} + \|\xd_k-\xdag\|\, \|\xd_{\kp}-\xdag\|
\Bigr)
\end{aligned}
\end{equation}
and
\begin{equation}\label{estBanach5_r1}
\begin{aligned}
&c \|T(\xd_{k+1}-\xdag)\| + \alpha_k \underline{c}\|\xd_{k+1}-\xdag\|^\pp\\
&\leq C \Bigl(\delta + \|\xd_k-\xdag\|^{3} + \|\xd_{\kp}-\xdag\|\,\|\xd_k-\xdag\|\\
&\quad+\Bigl(\|\xd_k-\xdag\|^{2}+\|\xd_k-\xdag\|+\|\xd_{\kp}-\xdag\|\Bigr)
\|\xd_{k+1}-\xdag\|
\Bigr)
\end{aligned}
\end{equation}
In case $\pp>1$, the elementary estimate \eqref{abeps} with $\pp$ instead of $\rr$ and $\epsilon=1$, 
$a=\left(\frac{\beta_k \underline{c}}{2C}\right)^{\frac{1}{\pp}}\|\xd_{\kp}-\xdag\|$, 
$b=\left(\frac{2C}{\beta_k \underline{c}}\right)^{\frac{1}{\pp}}\|\xd_k-\xdag\|$ implies 
\begin{equation}\label{estBanach6_r1}
\begin{aligned}
\|\xd_{\kp}-\xdag\|^\pp
\leq \frac{C(1,\pp)}{\underline{\beta}^\ppd} \left(\frac{2C}{\underline{c}}\right)^{\ppd}
\|\xd_k-\xdag\|^{\ppd}
+\frac{2C}{\underline{c}\underline{\beta}} \Bigl( \delta + \|\xd_k-\xdag\|^{2}\Bigr)
\end{aligned}
\end{equation}
and similarly 
\begin{equation}\label{estBanach7_r1}
\begin{aligned}
\|\xd_{k+1}-\xdag\|^\pp
\leq& \frac{C(1,\pp)}{\underline{\alpha}^\ppd} \left(\frac{2C}{\underline{c}}\right)^{\ppd}
\Bigl(\|\xd_k-\xdag\|^{2}+\|\xd_k-\xdag\|+\|\xd_{\kp}-\xdag\|\Bigr)^\ppd\\
&+\frac{2C}{\underline{c}\underline{\alpha}} \Bigl(\delta + \|\xd_k-\xdag\|^{3} + \|\xd_{\kp}-\xdag\|\,\|\xd_k-\xdag\|\Bigr)
\end{aligned}
\end{equation}
where we have used 
$\alpha_k\geq \underline{\alpha}$, $\beta_k\geq \underline{\beta}$.
If $\pp=1$, we have 
\begin{equation}\label{estBanach8_r1}
\begin{aligned}
\|\xd_{\kp}-\xdag\|
\leq \frac{C\Bigl( \delta + \|\xd_k-\xdag\|^{2}\Bigr)}{\underline{c}\underline{\beta}-C\|\xd_k-\xdag\|}
\end{aligned}
\end{equation}
and 
\begin{equation}\label{estBanach9_r1}
\begin{aligned}
\|\xd_{k+1}-\xdag\|
\leq &
\frac{C\Bigl(\delta + \|\xd_k-\xdag\|^{3} + \|\xd_{\kp}-\xdag\|\,\|\xd_k-\xdag\|\Bigr)}{\underline{c}\underline{\alpha}-C\Bigl(\|\xd_k-\xdag\|^{2}+\|\xd_k-\xdag\|+\|\xd_{\kp}-\xdag\|\Bigr)}\,.
\end{aligned}
\end{equation}
Inserting \eqref{estBanach6_r1} into \eqref{estBanach7_r1} and \eqref{estBanach8_r1} into \eqref{estBanach9_r1} we conclude that 
\[
\|\xd_{k+1}-\xdag\|\leq\mu_{k+1}\,,
\]
where $\mu_0=\|\xd_0-\xdag\|$,
\begin{equation}\label{muk}
\mu_{k+1}=\hat{C}\Bigl(\mu_k^\sigma +\delta^{\frac{1}{\pp}}\Bigr)\,,
\end{equation}
with $\hat{C}$ sufficiently large and 
\begin{equation}\label{sigma}
\sigma=\left\{\begin{array}{ll}
\frac{1}{\pp}\min\left\{2\ppd\,,\,\ppd\,,\,(\ppd)^2\,,\,3\,,\,1+\frac{\ppd}{\pp}\,,\,1+\frac{2}{\pp}\right\}
&\mbox{ if }\pp>1\\
3&\mbox{ if }\pp=1\end{array}\right\}=\frac{\pp+2}{\pp^2}\,,
\end{equation}
provided $\mu_k$ remains below some sufficiently small bound $\bar{\mu}>0$ (which will be guaranteed inductively by Lemma \ref{lem2} under smallness conditions on $\mu_0=\|x_0-\xdag\|$ and on $\delta$) such that in case $\pp=1$, 
$s\underline{c}-C\bar{\mu}>0$,
$\underline{c}-C\Bigl(\bar{\mu}^{2}+\bar{\mu}
+\frac{C}{s\underline{c}-C\bar{\mu}}\Bigl( \delta + \bar{\mu}^{2}\Bigr)\Bigr)>0$.
The requirement $\sigma>1$ resulting from the need of proving boundedness of $\mu_{k+1}$ according to \eqref{muk} with possibly large $\hat{C}$ translates to the condition $\pp<2$.
Now we make use of an elementary consequence of the recursion \eqref{muk}.
\begin{lemma}\label{lem2} 
For any $\hat{C}>0$, $\sigma>1$, $\pp\in[1,2)$, $\bar{\mu}\in(0,1]$, there exist $\bar{\mu}_0, \bar{\delta}>0$ such that for any $\delta\in [0,\bar{\delta}]$ and any $k_*\in\N$ we have the following:\\ 
Any sequence starting with $\mu_0\in[0,\bar{\mu}_0]$ and satisfying \eqref{muk} for all $k\in\{1,\ldots,k_*-1\}$ obeys the bound 
\[
\mu_{k+1}\leq 2^{-\sigma^{k+1}}
+C(\sigma) \delta^{\frac{1}{\pp}} \leq \bar{\mu} \mbox{ for all }k\leq k_*-1
\]
where $C(\sigma):=\sum_{m=0}^\infty 2^{-\sigma^m+1}$.
\end{lemma}
\begin{proof}
see the Appendix.
\end{proof}
So by setting $k=k_*-1$ according to \eqref{kstBanach_r1} we get
\[
\mu_{k_*} \leq (C(\sigma)+1) \delta^{\frac{1}{\pp}}\,,
\]
i.e., the stated convergence rate with noisy data.

If $\delta=0$ then \eqref{muk} directly provides us with convergence of $\mu_k$ to zero with convergence order $\sigma$.

\medskip

In the same manner the respective convergence result for the IRGNM in case $\rr=1$ can be seen: 
Namely, since the IRGNM corresponds to setting $x_{k+1}=x_{\kp}$, by \eqref{estBanach6_r1}, \eqref{estBanach8_r1} we have 
\[
\|
x_{k+1}^{\delta\mbox{\footnotesize IRGNM}}-\xdag\|\leq\tilde{\mu}_{k+1}\,,
\]
where $\tilde{\mu}_0=\|\xd_0-\xdag\|$,
\begin{equation*}
\tilde{\mu}_{k+1}=\tilde{\hat{C}}\Bigl(\tilde{\mu}_k^{\tilde{\sigma}} +\delta^{\frac{1}{\pp}}\Bigr)\,,
\end{equation*}
with $\tilde{\sigma}=\frac{1}{\pp}\min\{\ppd,2\}$.
The requirement $\tilde{\sigma}>1$ again translates to $ p<2 $, which entails that in fact $\tilde{\sigma}=\frac{2}{\pp}$ and Lemma \ref{lem2} yields the claimed result.

\section{Numerical experiments}\label{secNumTests}
We now show results of numerical tests with a Matlab implementation of method \eqref{HalleyIRGNMBanach} for the test example of identifying $c$ in
\[\begin{array}{rcll}
-\Delta u+\Upsilon(c)u&=&f&\mbox{ in } \Omega\\
u&=&g&\mbox{ on } \partial\Omega
\end{array}
\]
from measurements $y=Cu$ of $u$, 
where 
$\Upsilon(\lambda)=\frac12\lambda^2\chf _{[-\bar{c},\bar{c}]}
+\frac12\bar{c}(2|\lambda|-\bar{c})\chf _{\R\setminus[-\bar{c},\bar{c}]}$, so that $\Upsilon\in W^{2,\infty}(\R)$ and for the potential  $\Upsilon(c)$ nonnegativity and $\ppp$- integrability is guaranteed  if $c\in L^\ppp(\Omega)$: 
\[ 
\begin{aligned}
\|\Upsilon(c)\|_{L^\ppp}^\ppp=2^{-\ppp}\Bigl(
\int_{\{|c|\leq \bar{c}\}} |c|^{2\ppp}\, dx
+\int_{\{|c|> \bar{c}\}} |\bar{c}(2|c|-\bar{c})|^{\ppp}\, dx\Bigr)\\
\leq 2^{-\ppp}\Bigl(\bar{c}^\ppp
\int_{\{|c|\leq \bar{c}\}} |c|^{\ppp}\, dx
+(2\bar{c})^\ppp\int_{\{|c|> \bar{c}\}} |c|^{\ppp}\, dx\Bigr)
\leq  \bar{c}^{\ppp} \|c\|_{L^\ppp}^\ppp\,.
\end{aligned}
\]
For the forward operator $F=C\circ G$ with $G:L^\ppp(\Omega)\to W^{2,\ppp}(\Omega)$ and $C:W^{2,\ppp}(\Omega)\to Z$ some linear observation operator mapping in to some Banach space $Z$ we get $F'(c)h=CG'(c)h$, $F''(c)(h,l)=CG''(c)(h,l)$ with $v^1=G'(c)h$, $v^2=G''(c)(h,l)$ solving 
\begin{equation}\label{Fpccexample} 
\begin{array}{rcll}
-\Delta v^1+\Upsilon(c)v^1&=&-\Upsilon'(c)hG(c)&\mbox{ in } \Omega\\
v^1&=&0&\mbox{ on } \partial\Omega\,,
\end{array}
\end{equation}
\begin{equation}\label{Fppcexample} 
\begin{array}{rcll}
-\Delta v^2+\Upsilon(c)v^2&=&-\Upsilon'(c)hG'(c)l-\Upsilon'(c)lG'(c)h-\Upsilon''(c)hlG(c)&\mbox{ in } \Omega\\
v^2&=&0&\mbox{ on } \partial\Omega\,.
\end{array}
\end{equation}
Twice differentiability and Lipschitz continuity of $F''$ can be shown analogously to Example 3.1 in \cite{EKN89}, see also Example 2 in \cite{HalleyNM}.\\
Numerical tests are here done for
\begin{equation}\label{coscos}
c(x_1,x_2)=1+\frac52\xi(1-\cos(4\pi x_1))(1-\cos(4\pi x_2))\chf _{(0,\frac12)^2}\,,
\end{equation}
$f\equiv 4000$, $g\equiv 10$, and as starting value we use $c_0\equiv1$. 

\begin{figure}
\begin{center}
\includegraphics[width=0.32\textwidth]{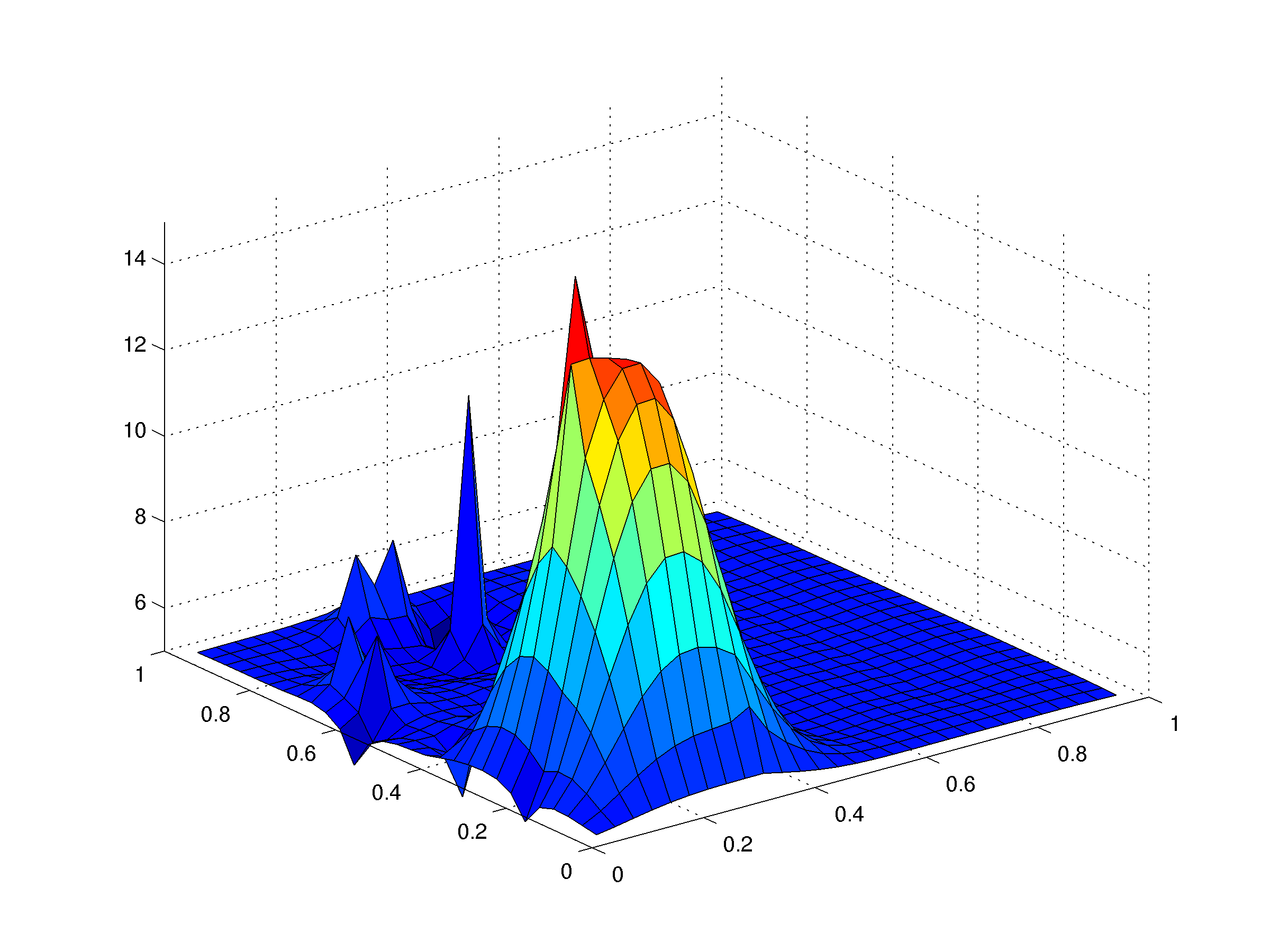}
\includegraphics[width=0.32\textwidth]{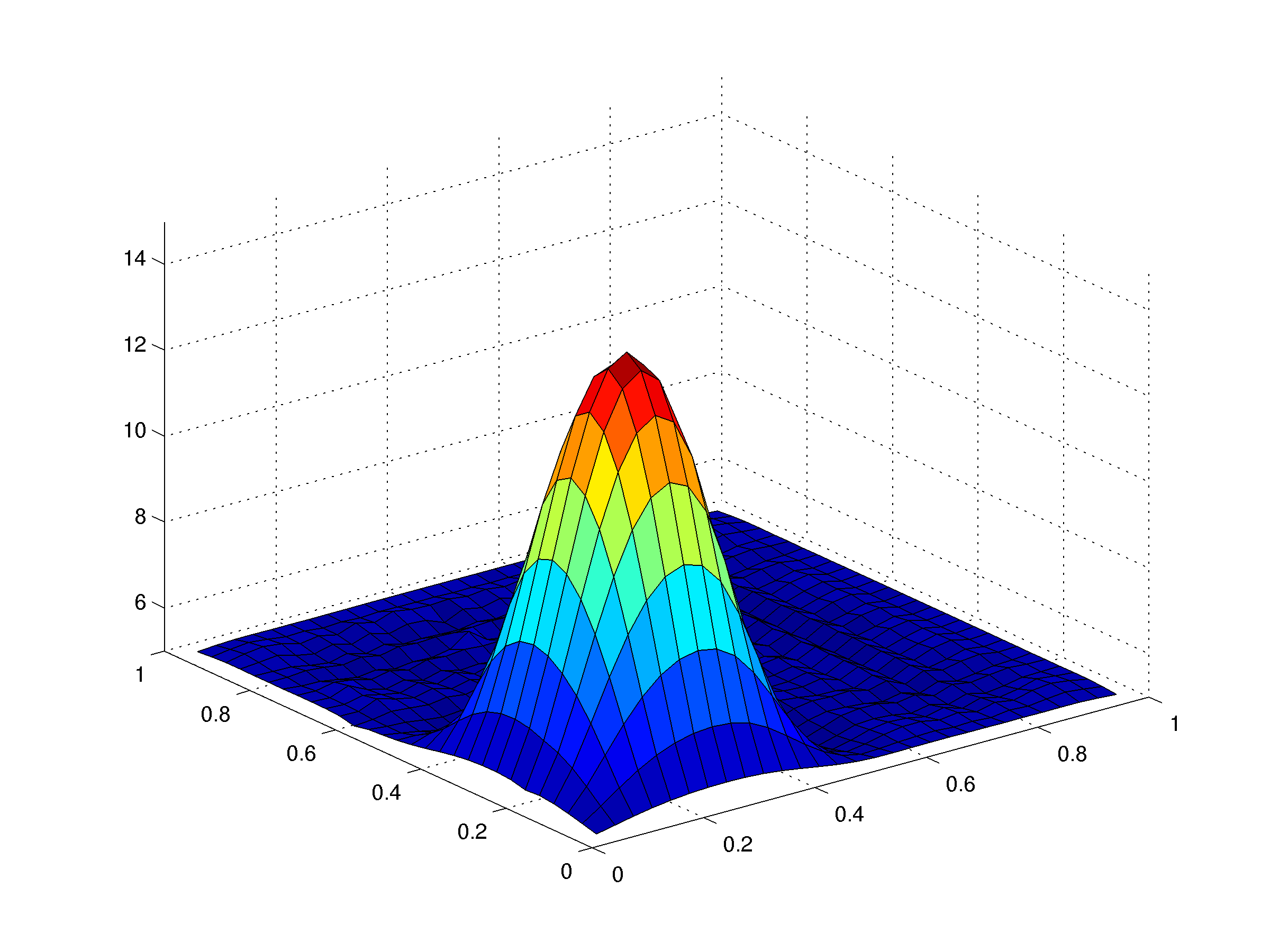}
\includegraphics[width=0.32\textwidth]{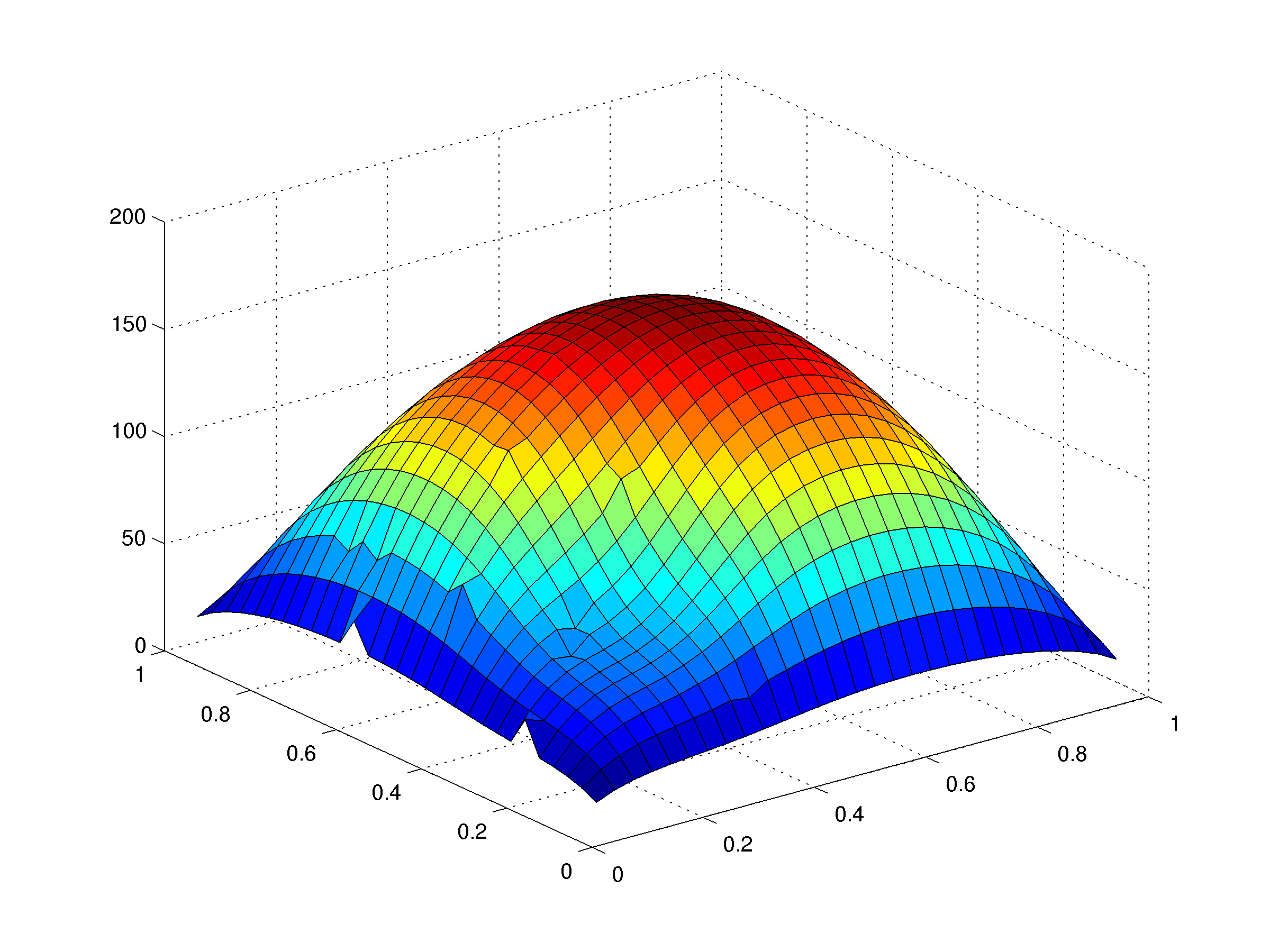}
\\
\includegraphics[width=0.32\textwidth]{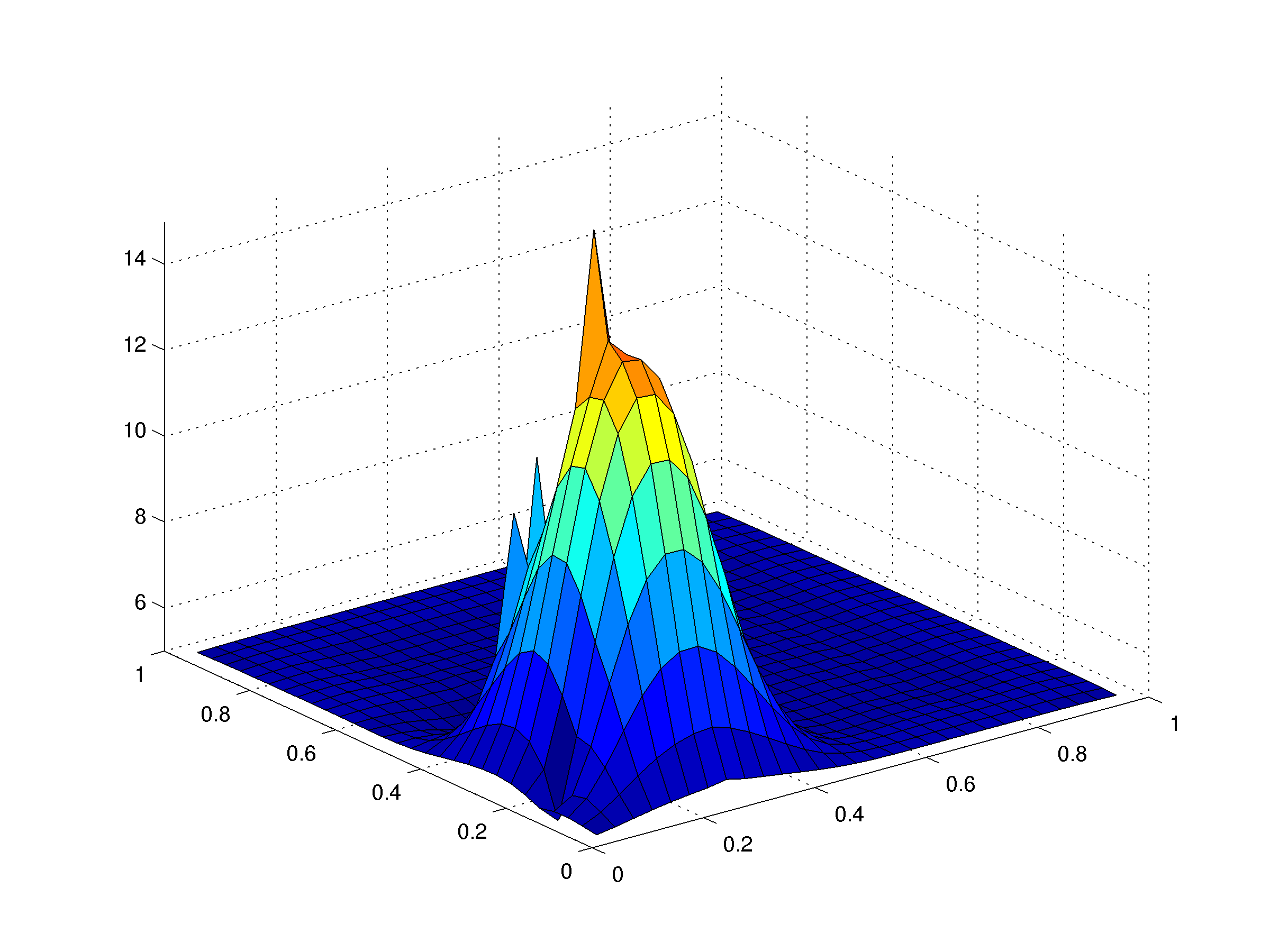}
\includegraphics[width=0.32\textwidth]{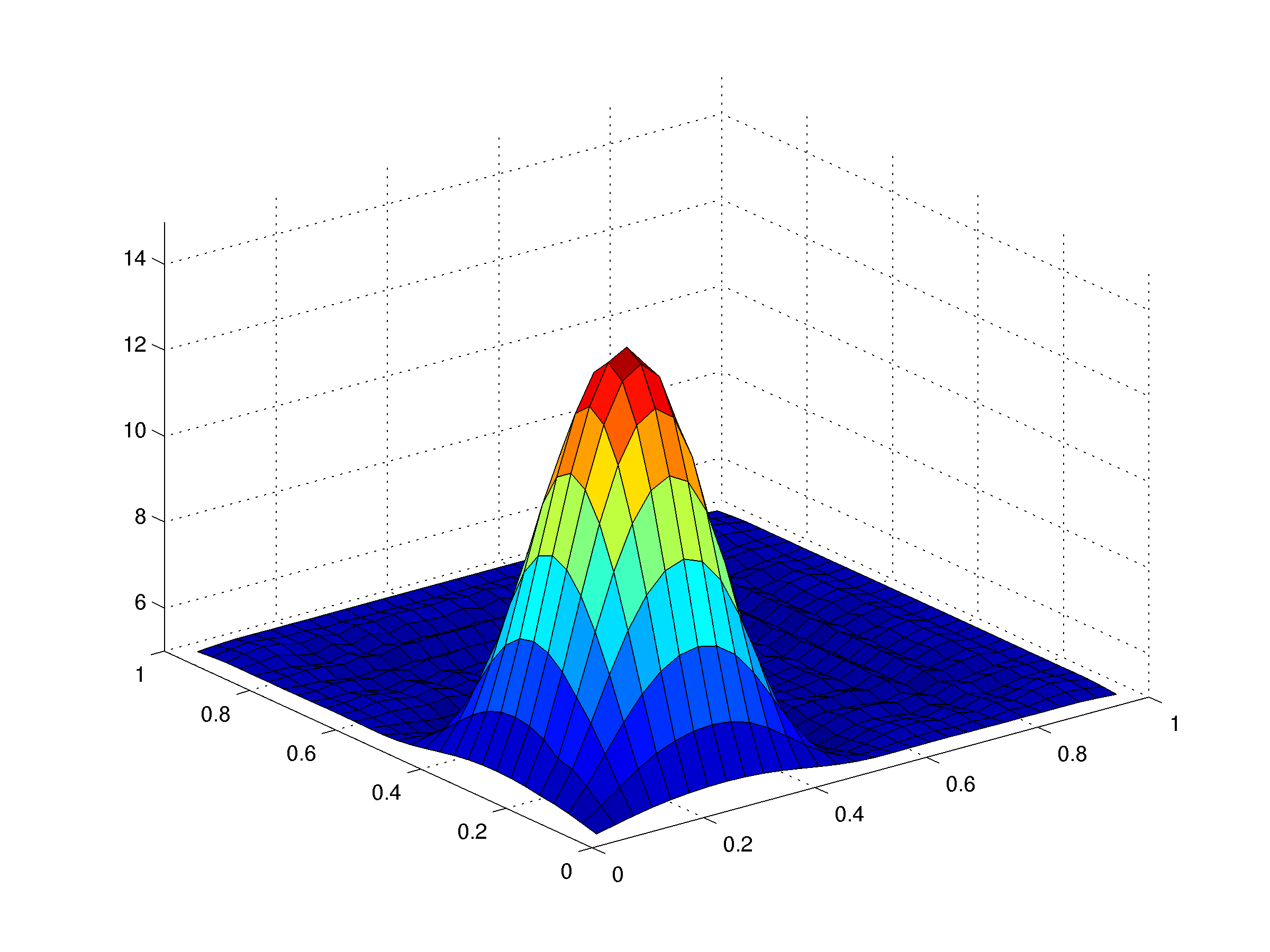}
\includegraphics[width=0.32\textwidth]{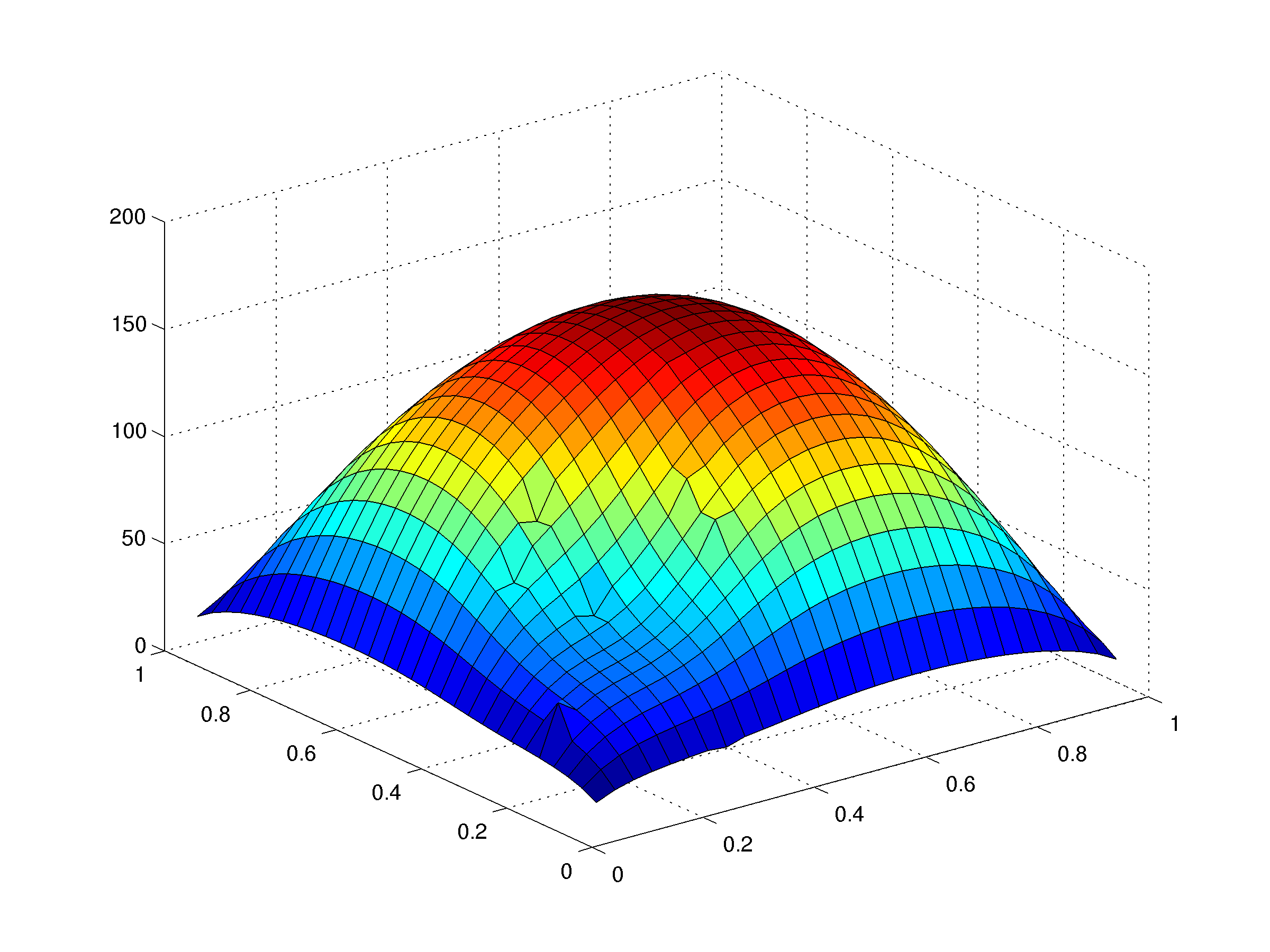}
\\
\includegraphics[width=0.32\textwidth]{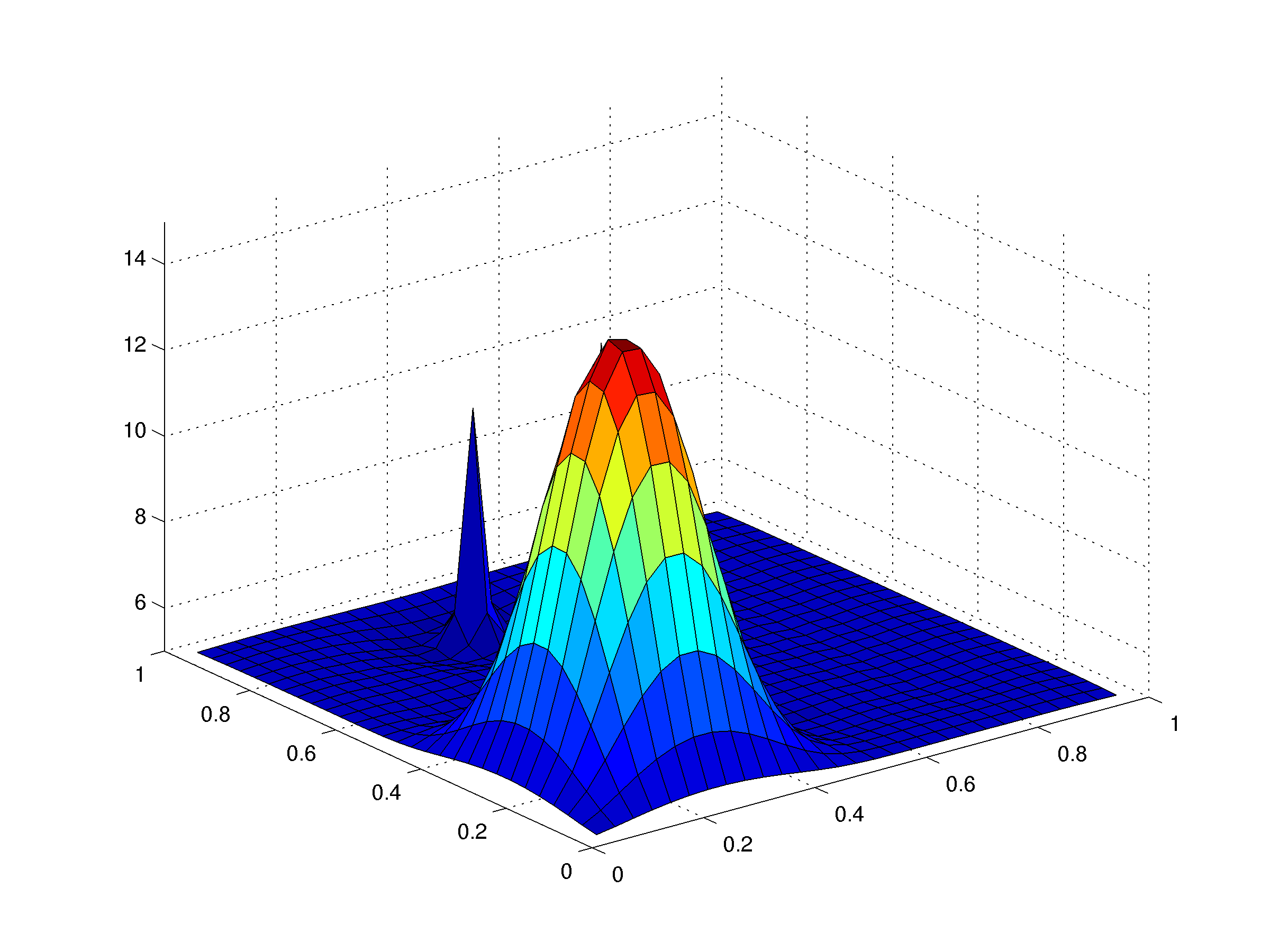}
\includegraphics[width=0.32\textwidth]{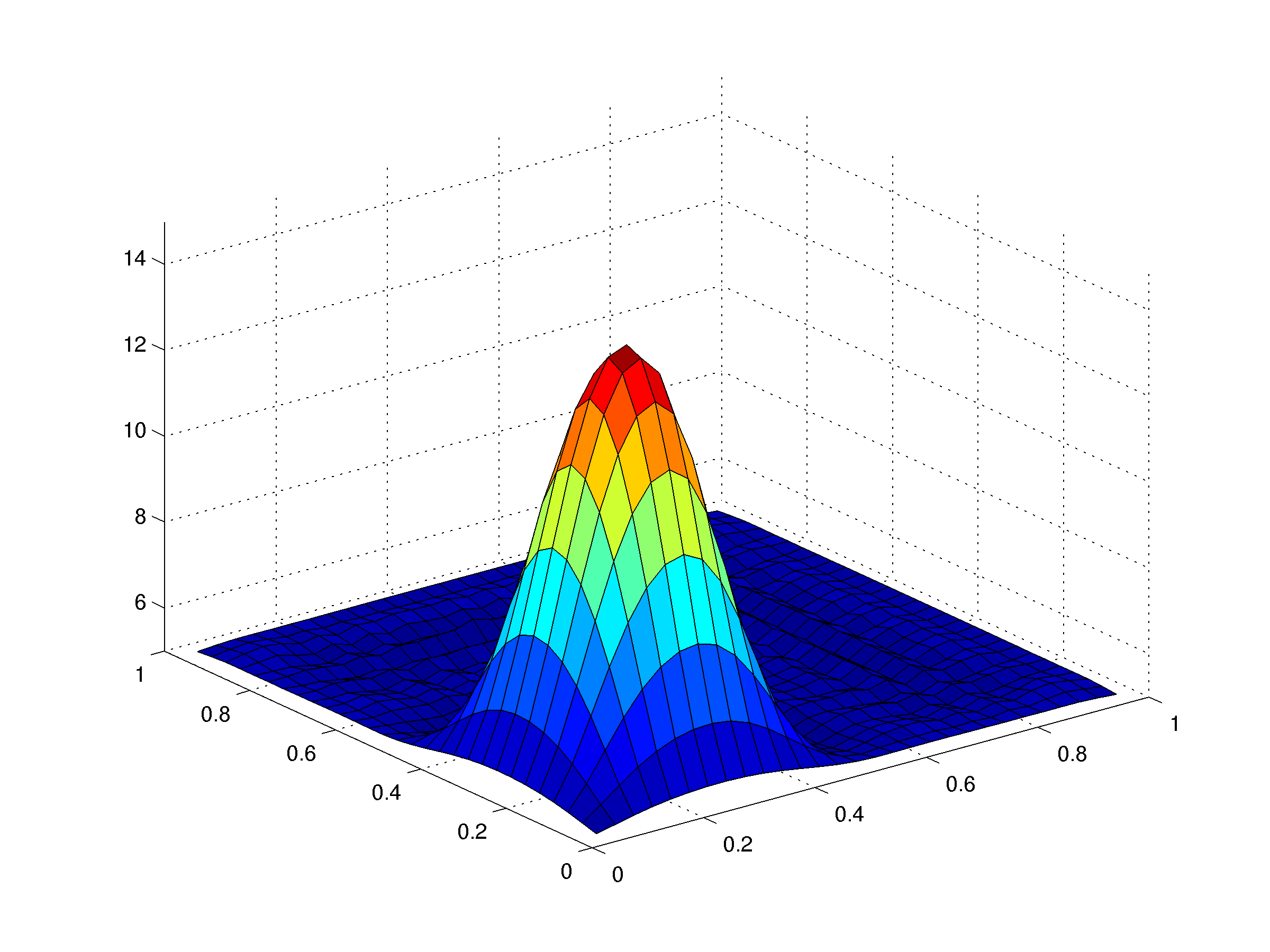}
\includegraphics[width=0.32\textwidth]{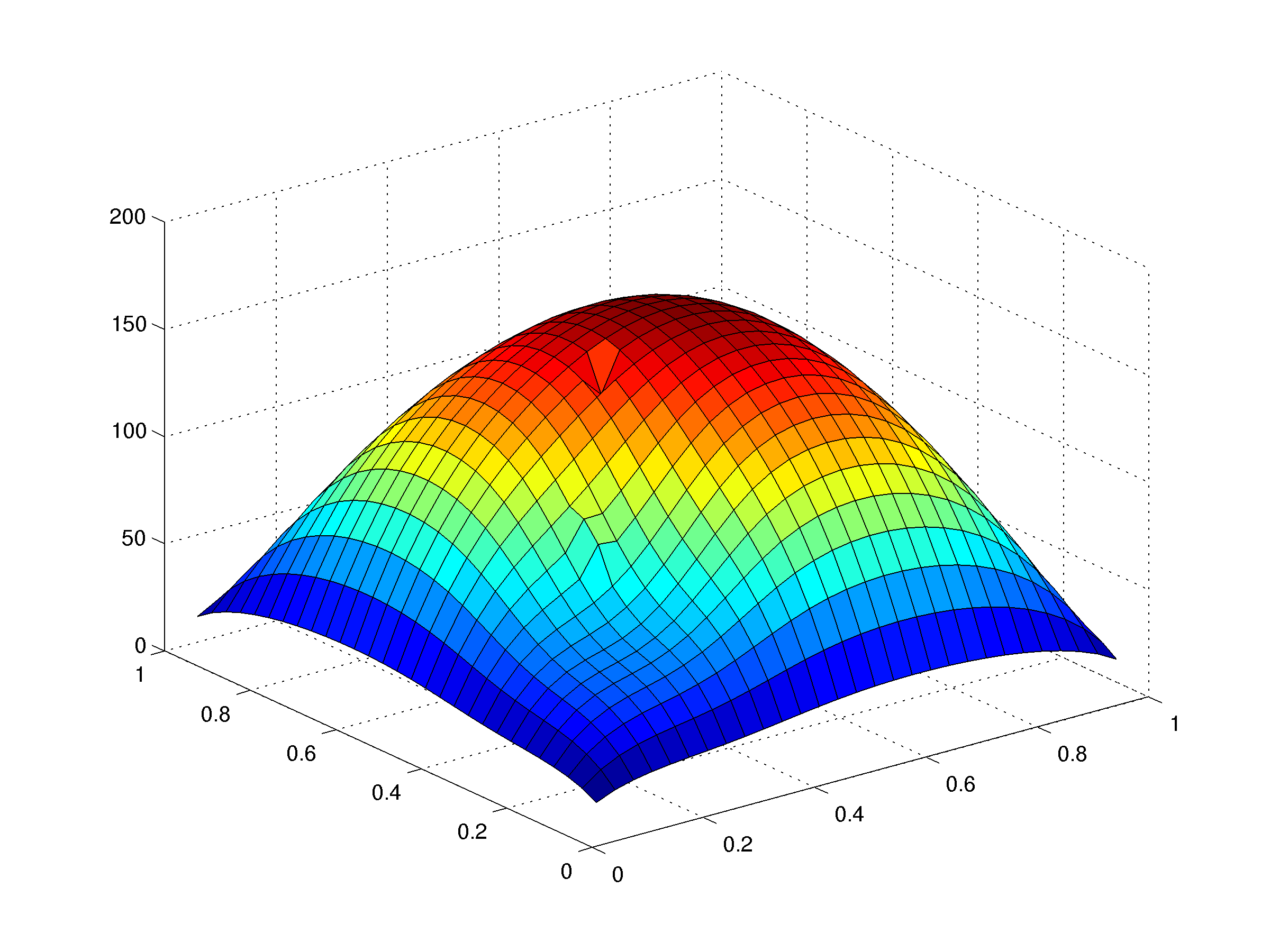}
\\
\includegraphics[width=0.32\textwidth]{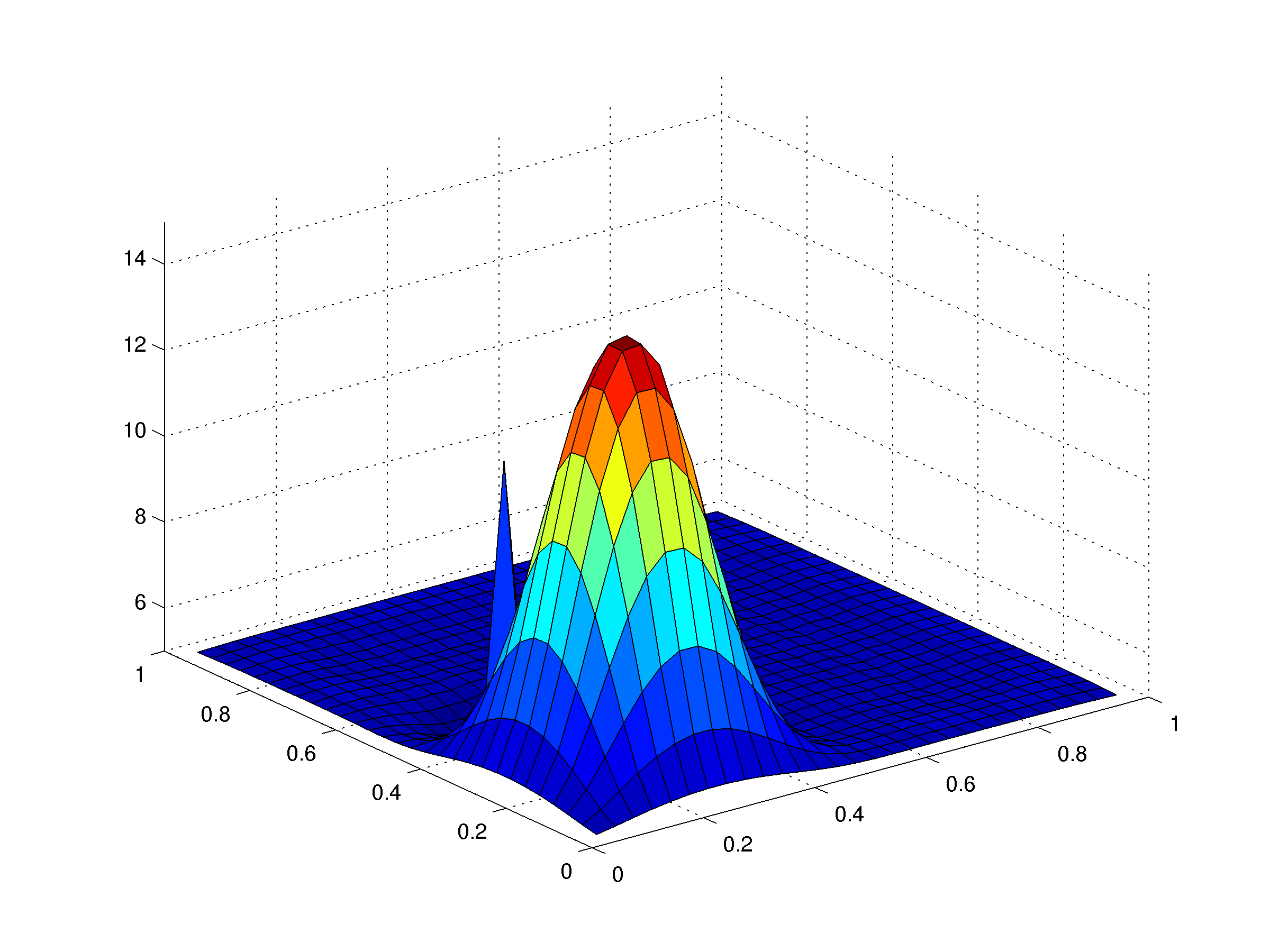}
\includegraphics[width=0.32\textwidth]{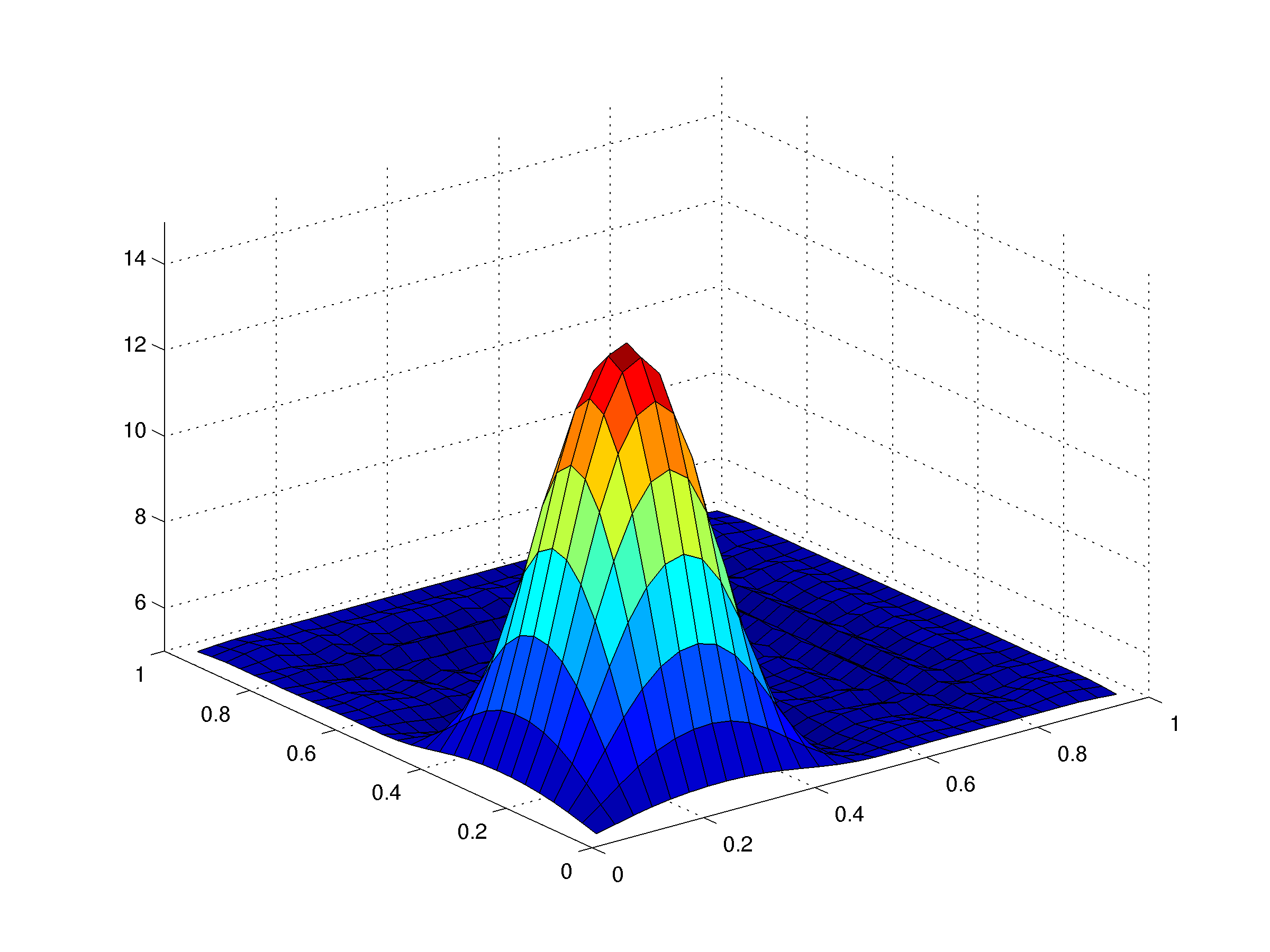}
\includegraphics[width=0.32\textwidth]{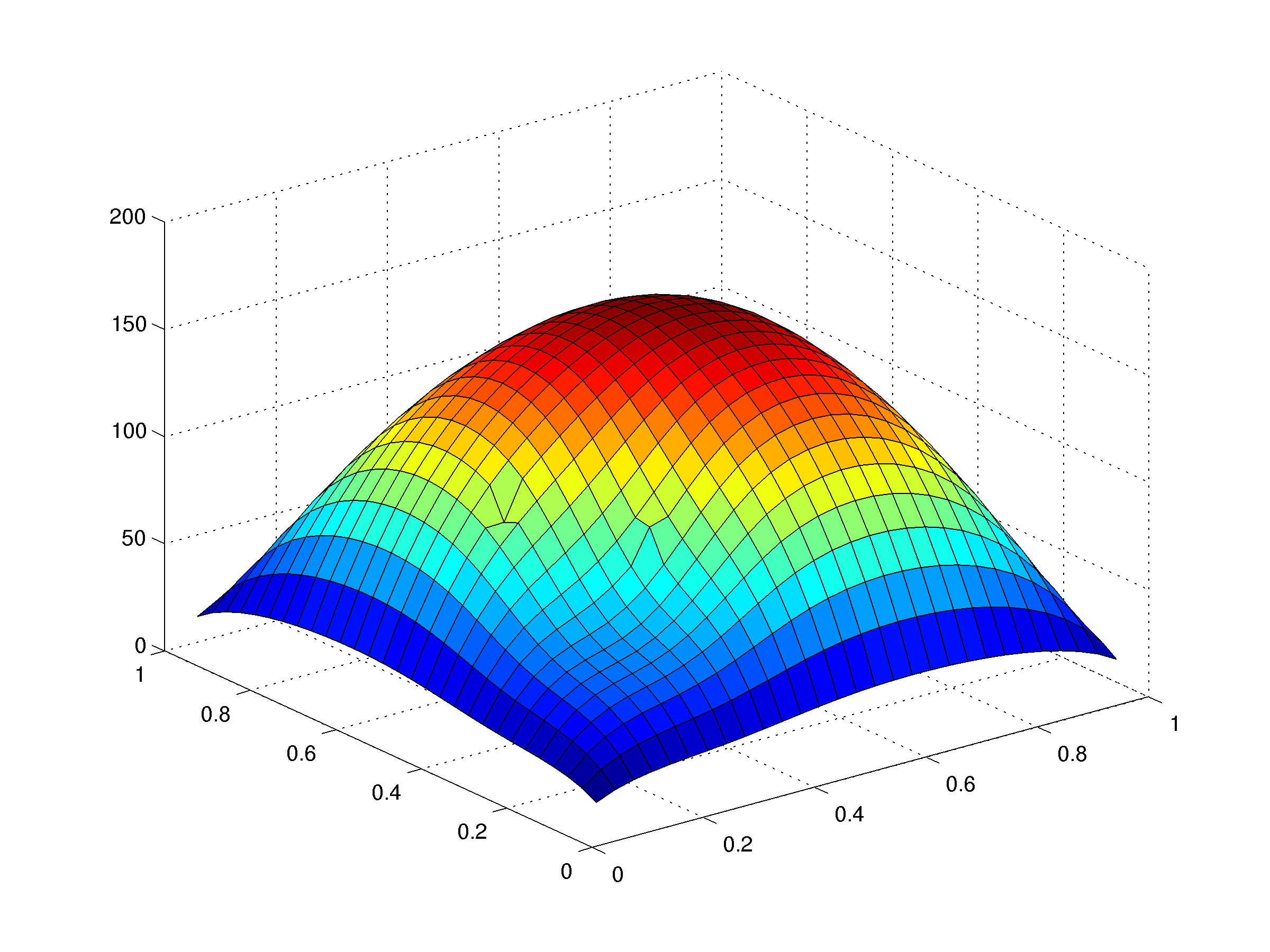}
\\
\includegraphics[width=0.32\textwidth]{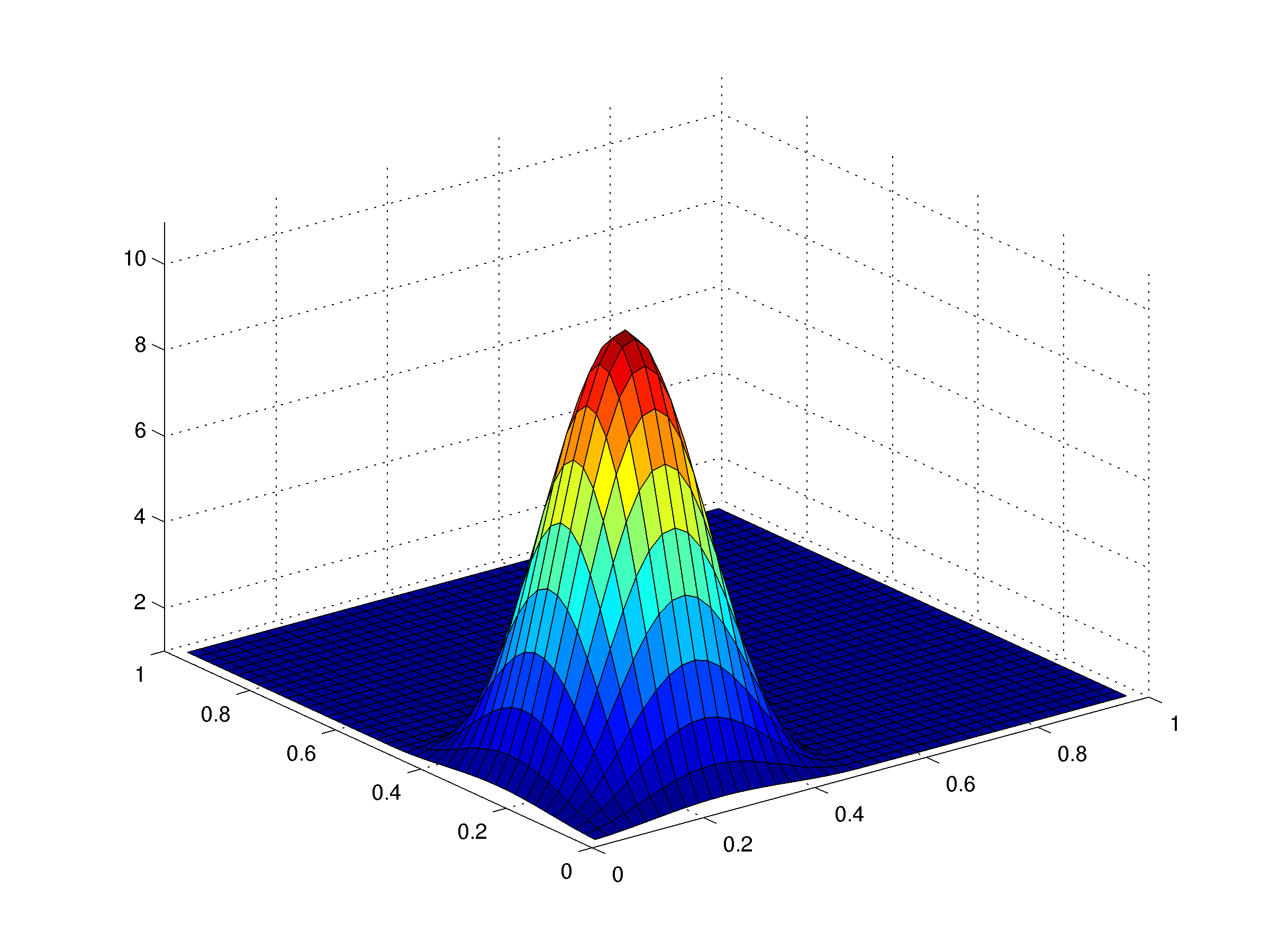}
\includegraphics[width=0.32\textwidth]{qex_coscos.png}
\includegraphics[width=0.32\textwidth]{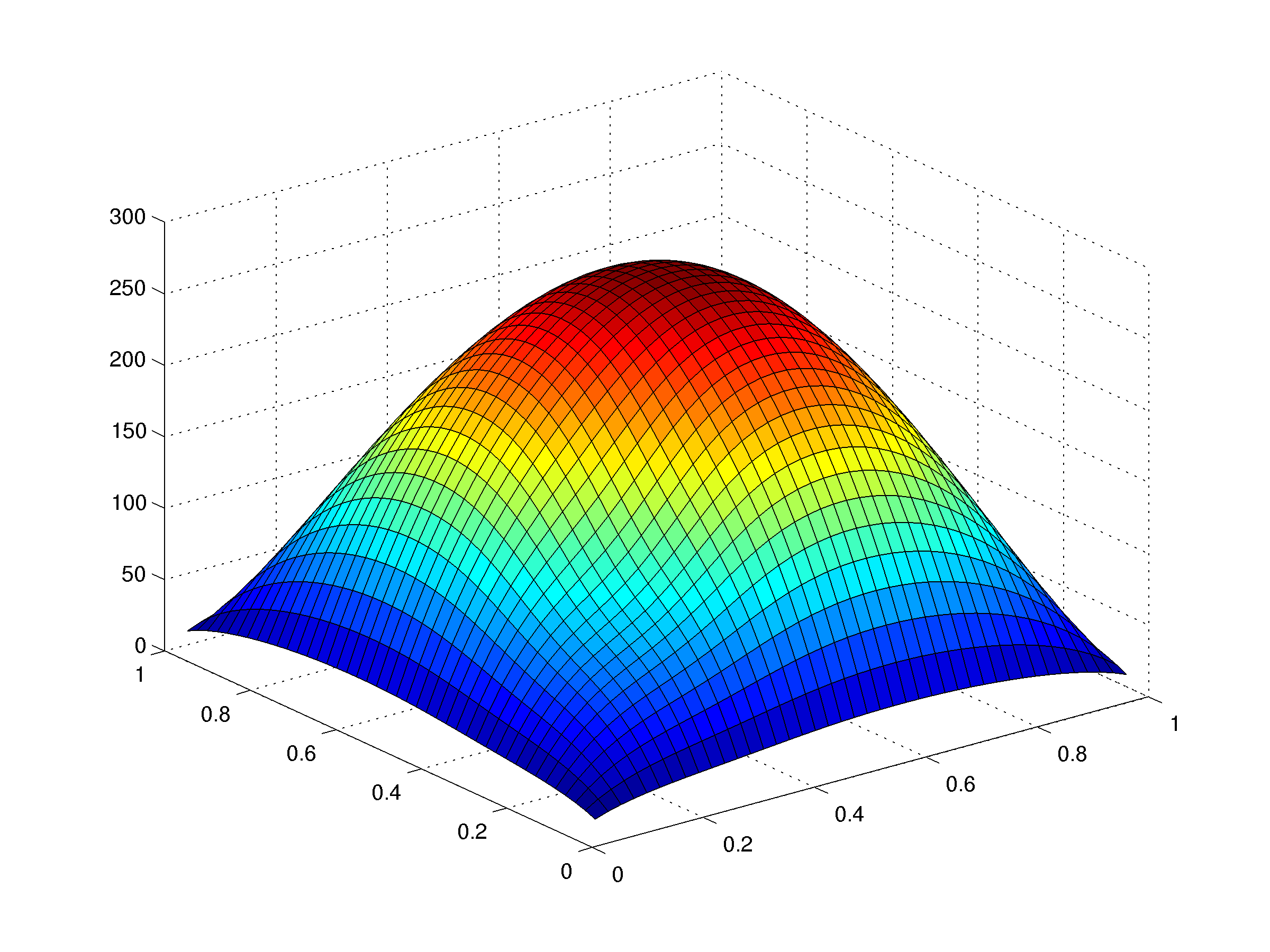}
\caption{Reconstructions with $Y=L^2$ (left) and $Y=L^{1.1}$ (middle) from data (right) with outliers noise of decreasing amount from top to bottom, compared to exact $c$ and $u$ (bottom row)
\label{fig_saltandpepper}
}
\end{center}
\end{figure}

A comparison with the IRGNM for this example has been carried out in the Hilbert space setting of \cite{HalleyNM}.
Here we consider non-Gaussian noise and compare performance of the Hilbert space version of Halley's method with the formulation in appropriate Banach spaces. It is well-known that $L^1$ 
data misfit terms are better suited than the $L^2$ norm in case of impulsive noise.
This can be seen also here in our tests with $Y=L^{1.1}$, $\rr=1.1$ 
%
as compared to  $Y=L^{2}$, $\rr=2$ (Figure \ref{fig_saltandpepper}).
Here the noise was generated by randomly (uniformly distributed) picking measurement points and perturbing their values by an amount of ten per cent of the maximal measurement value.  
Note that the 
ideal choice $Y=L^1$ 
is admissible by our theory but
would make the subproblems in each Halley step nonsmooth, which would require more sophisticated numerical techniques, (see, e.g., \cite{Clason})
than what we have implemented for our tests.
Using $L^{1+\epsilon}$ with $\epsilon>0$ may be viewed as a smooth approximation to the computationally hard case of $L^1$ or the space of Radon measures, cf. \cite{StrehlowKazimierski14}.

\section{Conclusions and remarks}\label{secConcl}
In this paper we have extended the IRGNM-Halley method from \cite{HalleyNM} to the general Banach space setting with possibly nonquadratic penalties and proven convergence rates under a particular source condition and with a
priori regularization parameter choice for a variety of exponents in the data misfit and regularization terms.

More general convergence rates, including convergence without rates, have yet to be shown. Such results might be obtained using approximate or variational source conditions. As soon as \eqref{scBanach} is violated, certainly stronger structural assumptions on $F$ will be needed to still establish convergence. It is not yet clear, though, how such conditions should be formulated to enable convergence proofs and still be satisfied for relevant applications.
In \cite{HettlichRundell} a tangential cone type condition was successfully used for proving convergence without source conditions. However, for the Levenberg-Marquardt type approach taken there, a monotonicity argument can be used, which does not apply to the IRGNM-type version considered here.
In the Hilbert space setting of \cite{HalleyNM}, we have proven convergence without (or with weaker) source conditions under a range invariance condition on $F'$, $F''$, (which is in some sense dual to the tangential cone condition). However, it is not yet clear how to carry out proofs under such conditions in a Banach space setting, where the classical Hilbert space spectral calculus is not available. Even for the first order version, i.e., the original IRGNM, this has yet to be done in non-Hilbert spaces.

Further research will therefore be concerned with providing answers to these open questions. 

\section{Appendix}
\begin{proof} (Lemma \ref{lem1})\\
For any fixed tentative constant $\bar{\gamma}>0$ 
we ask for existence of a $\bar{\Gamma}>0$ such that the implication 
\[
\Gamma^\pp\leq \phi(\bar{\gamma},\Gamma) \ \Rightarrow \ \Gamma\leq \bar{\Gamma}
\]
holds. By contraposition this is equivalent to existence of a $\bar{\Gamma}>0$ such that
\[
\Gamma> \bar{\Gamma} \ \Rightarrow \ \Gamma^\pp - \phi(\bar{\gamma},\Gamma) >0
\]
holds, which by continuity of $\phi$ implies 
\[
\lim_{\Gamma\to\infty} \Gamma^\pp - \phi(\bar{\gamma},\Gamma) >0\,.
\]
By inspection of the function 
$\phi(\bar{\gamma},\cdot):\Gamma\mapsto a+b\bar{\gamma}^{2\rr}+c\bar{\gamma}^{\rr}\Gamma^{\rr}$,
the latter can be easily seen to be equivalent to 
\begin{equation}\label{pger}
\pp\geq\rr
\end{equation}
(and additionally $c \bar{\gamma}^\rr \leq 1$ in case $\pp=\rr$). After having derived a necessary condition relating $\pp$ and $\rr$ we now return to the inequality
\begin{equation} \label{estBanach9}
\Gamma^\pp\leq \phi(\bar{\gamma}, \Gamma)
\end{equation}
and compute a resulting explicit upper estimate of $\Gamma$ in terms of $\bar{\gamma}$ by distinction between the cases $\Gamma\leq 1$ and $\Gamma> 1$, the latter by \eqref{estBanach9} resulting in 
\[
a+b\bar{\gamma}^{2\rr}\geq \psi(\Gamma^\pp)
=\psi(1)+\psi'(1+\theta(\Gamma^\pp-1)) (\Gamma^\pp-1)
\]
with 
\[
\psi(\lambda):=\lambda-c\bar{\gamma}^\rr\lambda^{\frac{\rr}{\pp}}\,, \quad
\psi'(\lambda)= 1-c\frac{\rr}{\pp}\frac{\bar{\gamma}^\rr}{\lambda^{(1-\frac{\rr}{\pp})}}
\geq 1-c\frac{\rr}{\pp}\bar{\gamma}^\rr \mbox{ for }\lambda\geq 1
\]
(where we have used \eqref{pger}), hence 
\[
\Gamma^\pp\leq 1+ \frac{a+b\bar{\gamma}^{2\rr}-1+c\bar{\gamma}^\rr}{1-c\frac{\rr}{\pp}\bar{\gamma}^\rr}\,,
\]
provided $c\frac{\rr}{\pp}\bar{\gamma}^\rr<1$, which altogether gives 
\[
\Gamma\leq 
\left(1+ \frac{a+b\bar{\gamma}^{2\rr}+c\bar{\gamma}^\rr}{1-c\frac{\rr}{\pp}\bar{\gamma}^\rr}\right)^{\frac{1}{\pp}}
\]
in either of the two cases $\Gamma\leq/>1$.
Inserting this into 
\[
\gamma^\pp\leq \Phi(\bar{\gamma},\gamma,\Gamma)
\]
 yields
\begin{equation}\label{estBanach10}
A(\bar{\gamma})\geq \Psi(\gamma^\pp)
=\Psi(\lambda_0)+\Psi'(\lambda_0+\theta(\Gamma^\pp-\lambda_0)) (\gamma^\pp-\lambda_0)
\end{equation}
for some $\theta\in[0,1]$ with 
\[
A(\bar{\gamma}):=
d+ e \bar{\gamma}^{3\rr}+ f \bar{\gamma}^\rr
\Bigl(1+ \frac{a+b\bar{\gamma}^{2\rr}+c\bar{\gamma}^\rr}{1-c\frac{\rr}{\pp}\bar{\gamma}^\rr}\Bigr)^\frac{\rr}{\pp}
\]
\[
B(\bar{\gamma}):=h \bar{\gamma}^{2\rr}
+ i \bar{\gamma}^{\rr}
+ j \Bigl(1+ \frac{a+b\bar{\gamma}^{2\rr}+c\bar{\gamma}^\rr}{1-c\frac{\rr}{\pp}\bar{\gamma}^\rr}\Bigr)^{\frac{\rr}{\pp}}
\]
\[
\Psi(\lambda):=\lambda-B(\bar{\gamma}) \lambda^{\frac{\rr}{\pp}}\,, \quad
\Psi'(\lambda)=1-\frac{\rr}{\pp}\frac{B(\bar{\gamma})}{\lambda^{1-\frac{\rr}{\pp}}}
\]
Thus, similarly to above, by distinction between the cases $\gamma^\pp</\ge \lambda_0$ we can estimate
\[
\gamma^\pp<\lambda_0\mbox{ or }
\lambda_0\leq\gamma^\pp\leq\lambda_0+\frac{A(\bar{\gamma})-\lambda_0+B(\bar{\gamma})\lambda_0^{\frac{\rr}{\pp}}}{1-\frac{\rr}{\pp}\frac{B(\bar{\gamma})}{\lambda_0^{1-\frac{\rr}{\pp}}}}
\]
It remains to show that the right hand side of this inequality can be bounded by $\bar{\gamma}^\pp$, using a proper choice of $\bar{\gamma}>0$ and $\lambda_0>0$. We do so by setting $\lambda_0=\left(\frac{\bar{\gamma}}{3}\right)^\pp$, so that it remains to show that 
\[
\begin{aligned}
&\gamma^\pp<\left(\frac{\bar{\gamma}}{3}\right)^\pp\mbox{ or }\\
&0\leq\gamma^\pp-\left(\frac{\bar{\gamma}}{3}\right)^\pp
\leq \frac{A(\bar{\gamma})-\left(\frac{\bar{\gamma}}{3}\right)^\pp+B(\bar{\gamma})(\left(\frac{\bar{\gamma}}{3}\right)^\pp)^{\frac{\rr}{\pp}}}{1-\frac{\rr}{\pp}\frac{B(\bar{\gamma})}{(\left(\frac{\bar{\gamma}}{3}\right)^\pp)^{1-\frac{\rr}{\pp}}}}\leq (3^\pp-1)\left(\frac{\bar{\gamma}}{3}\right)^\pp\,,
\end{aligned}
\]
i.e., unless $\gamma^\pp\leq \left(\frac{\bar{\gamma}}{3}\right)^\pp$ happens to hold (in which case we would already be finished)
\[
0\leq A(\bar{\gamma})-\left(\frac{\bar{\gamma}}{3}\right)^\pp
+B(\bar{\gamma})\left(\frac{\bar{\gamma}}{3}\right)^{\rr} 
\leq (3^\pp-1)\left(\left(\frac{\bar{\gamma}}{3}\right)^\pp-\frac{\rr}{\pp} B(\bar{\gamma})\left(\frac{\bar{\gamma}}{3}\right)^\rr\right)
\]
must be shown.
Considering the asymptotic behavior as $\bar{\gamma}\to0$ yields the requirement
\begin{equation}\label{asymp}
\begin{aligned}
0\leq& d+e\bar{\gamma}^{3\rr}+f\bar{\gamma}^{\rr}
\left(1+a\pm O(\bar{\gamma}^{\rr})\right)^{\frac{\rr}{\pp}}
-\left(\frac{\bar{\gamma}}{3}\right)^\pp\\
&+ \left( h\bar{\gamma}^{2\rr}+i\bar{\gamma}^{\rr} + j
\left(1+a\pm O(\bar{\gamma}^{\rr})\right)^{\frac{\rr}{\pp}}\right)
\left(\frac{\bar{\gamma}}{3}\right)^\rr\\
\leq& (3^\pp-1) \left(\left(\frac{\bar{\gamma}}{3}\right)^\pp-\frac{\rr}{\pp} 
\left( h\bar{\gamma}^{2\rr}+i\bar{\gamma}^{\rr} + j
\left(1+a+ O(\bar{\gamma}^{\rr})\right)^{\frac{\rr}{\pp}}\right)
\left(\frac{\bar{\gamma}}{3}\right)^\rr\right)\,.
\end{aligned}
\end{equation}
This shows that we have to decrease $d,e,f,h,j$ depending on $\bar{\gamma}$, i.e., we assume that we can choose 
\[
d=\frac32 \left(\frac{\bar{\gamma}}{3}\right)^\pp\,, \quad
e=o(\bar{\gamma}^{\pp-3\rr})\,, \quad
f=o(\bar{\gamma}^{\pp-\rr})\,, \quad
h=o(\bar{\gamma}^{\pp-3\rr})\,, \quad
j=o(\bar{\gamma}^{\pp-\rr})\,.
\] 
Also for the $i$ term we need $i\bar{\gamma}^{2\rr}=o(\bar{\gamma}^{\pp})$, which can be achieved by assuming 
\[ \pp<2\rr \]
(note that by $\rr\geq1$ this is less restrictive than assuming $\pp<2$ in order to make $i=\frac{C}{\underline{c}}q^{\frac{\rrd}{\ppd}} \alpha_0^\frac{\rrd(2-\pp)}{\pp}$ small).
These choices render \eqref{asymp} an asymptotic estimate of the form
\[
0\leq \frac12 \left(\frac{\bar{\gamma}}{3}\right)^\pp\pm o(\bar{\gamma}^\pp)
\leq (3^\pp-1) \left(\frac{\bar{\gamma}}{3}\right)^\pp -o(\bar{\gamma}^\pp)
\]
which is obviously feasible, so that the desired estimate
\[\gamma\leq \bar{\gamma}\]
can be achieved by choosing $\bar{\gamma}$ sufficiently small.
\end{proof}

\medskip

\begin{proof} (Lemma \ref{lem2})\\
For any $l\leq k$, the estimate
\begin{equation}\label{estmukl}
\begin{aligned}
\mu_{k+1}\leq&\hat{C}^{\frac{\sigma^{l+1}-1}{\sigma-1}} 2^{\frac{\sigma}{\sigma-1}(\sigma^l-1)-l}
\mu_{k-l}^{\sigma^{l+1}}\\
&+\sum_{m=0}^l \hat{C}^{\frac{\sigma^{m+1}-1}{\sigma-1}} 2^{\frac{\sigma}{\sigma-1}(\sigma^m-1)-m}
\delta^{\frac{\sigma^m}{\pp}}
\end{aligned}
\end{equation}
which can be seen by induction and the elementary estimate 
$\left( a+b \right)^\lambda\leq 2^{\lambda-1}a^\lambda+2^{\lambda-1}b^\lambda$ for $a,b\geq0$, $\lambda\geq1$. Namely, from \eqref{muk} we have 
\[
\begin{aligned}
\mu_{k+1}=&\hat{C}\Bigl( \Bigl(\hat{C}\Bigl(\mu_{k-1}^\sigma +\delta^{\frac{1}{\pp}}\Bigr)\Bigr)^\sigma +\delta^{\frac{1}{\pp}}\Bigr)\\
\leq& \hat{C} 2^{\sigma-1} \hat{C}^\sigma \mu_{k-1}^{\sigma^2}
+\hat{C} 2^{\sigma-1} \hat{C}^\sigma \delta^{\frac{\sigma}{\pp}}
+\hat{C} \delta^{\frac{1}{\pp}}
\end{aligned}
\]
which is just \eqref{estmukl} with $l=1$.
To carry out the induction step we again use \eqref{muk} with $k$ replaced by $k-l-1$ in \eqref{estmukl} to obtain
\[
\begin{aligned}
\mu_{k+1}\leq&\hat{C}^{\frac{\sigma^{l+1}-1}{\sigma-1}} 2^{\frac{\sigma}{\sigma-1}(\sigma^l-1)-l}
\Bigl(\hat{C}\Bigl(\mu_{k-l-1}^\sigma +\delta^{\frac{1}{\pp}}\Bigr)\Bigr)^{\sigma^{l+1}}\\
&+\sum_{m=0}^l \hat{C}^{\frac{\sigma^{m+1}-1}{\sigma-1}} 2^{\frac{\sigma}{\sigma-1}(\sigma^m-1)-m}
\delta^{\frac{\sigma^m}{\pp}}\\
\leq&\hat{C}^{\frac{\sigma^{l+1}-1}{\sigma-1}} 2^{\frac{\sigma}{\sigma-1}(\sigma^l-1)-l}
2^{\sigma^{l+1}-1}\hat{C}^{\sigma^{l+1}}\mu_{k-l-1}^{\sigma^{l+2}} \\
&+\hat{C}^{\frac{\sigma^{l+1}-1}{\sigma-1}} 2^{\frac{\sigma}{\sigma-1}(\sigma^l-1)-l}
2^{\sigma^{l+1}-1}\hat{C}^{\sigma^{l+1}} \delta^{\frac{\sigma^{l+1}}{\pp}}\\
&+\sum_{m=0}^l \hat{C}^{\frac{\sigma^{m+1}-1}{\sigma-1}} 2^{\frac{\sigma}{\sigma-1}(\sigma^m-1)-m}
\delta^{\frac{\sigma^m}{\pp}}\\
=&\hat{C}^{\frac{\sigma^{l+2}-1}{\sigma-1}} 2^{\frac{\sigma}{\sigma-1}(\sigma^{l+1}-1)-l-1}
\mu_{k-l-1}^{\sigma^{l+2}}\\
&+\sum_{m=0}^{l+1} \hat{C}^{\frac{\sigma^{m+1}-1}{\sigma-1}} 2^{\frac{\sigma}{\sigma-1}(\sigma^m-1)-m}
\delta^{\frac{\sigma^m}{\pp}}
\end{aligned}
\]
which completes the proof of \eqref{estmukl}. We now use $l=k$ in \eqref{estmukl} to conclude that 
\[
\begin{aligned}
\mu_{k+1}\leq&\hat{C}^{\frac{\sigma^{k+1}-1}{\sigma-1}} 2^{\frac{\sigma}{\sigma-1}(\sigma^k-1)-k}
\mu_{0}^{\sigma^{k+1}}\\
&+\sum_{m=0}^k \hat{C}^{\frac{\sigma^{m+1}-1}{\sigma-1}} 2^{\frac{\sigma}{\sigma-1}(\sigma^m-1)-m}
\delta^{\frac{\sigma^m}{\pp}}\\
\leq& \left((2\hat{C})^{\frac{1}{\sigma-1}}\mu_0\right)^{\sigma^{k+1}}
+\left(\sum_{m=0}^k \left((2\hat{C})^{\frac{\sigma}{\sigma-1}}
\delta^{\frac{1}{\pp}}\right)^{\sigma^m-1}\right)
\delta^{\frac{1}{\pp}}\\
\leq& 2^{-\sigma^{k+1}} \bar{\mu}
+\left(\sum_{m=0}^k 2^{-\sigma^m+1}\right)
\delta^{\frac{1}{\pp}} 
\leq \bar{\mu}
\leq 1
\end{aligned}
\]
under the smallness assumptions
\begin{equation} \label{smallnessBanach}
\mu_0\leq \frac{\bar{\mu}}{2} (2\hat{C})^{-\frac{1}{\sigma-1}}\,, \quad
\delta\leq\bar{\delta}:=\min\{ \bar{\mu}\frac{1-2^{-\sigma^2}}{C(\sigma)}\,, \ \frac12(2\hat{C})^{-\frac{\pp\sigma}{\sigma-1}}\}\,,
\end{equation}
where $C(\sigma):=\sum_{m=0}^\infty 2^{-\sigma^m+1}$.
\end{proof}

\section*{Acknowledgment}
Financial support by the Austrian Science Fund FWF under grant P24970 is gratefully acknowledged.

\end{document}